\documentclass[11pt,final,3p]{elsarticle}
\usepackage{amssymb}
\usepackage{amsthm}
\usepackage{amsmath}
\usepackage{latexsym}
\usepackage{amsmath}
\usepackage{xcolor}
\usepackage{graphicx}
\usepackage{indentfirst}
\usepackage[OT2,OT1]{fontenc}
\usepackage{mathrsfs}
\biboptions{numbers,sort&compress}
\usepackage{pdfsync}
\usepackage{hyperref}
\hypersetup{hypertex=true,
	colorlinks=true,
	linkcolor=blue,
	anchorcolor=blue,
	citecolor=blue}
\allowdisplaybreaks
\usepackage[utf8]{inputenc}

\renewcommand{\thefootnote}{\arabic{footnote}}

\numberwithin{equation}{section}

\newcommand\blfootnote[1]{%
	\begingroup
	\renewcommand\thefootnote{}\footnote{#1}%
	\addtocounter{footnote}{-1}%
	\endgroup
}
\journal{ }
\begin{document}
\begin{frontmatter}
\title{ McKean-Vlasov SPDEs with coefficients exhibiting locally weak monotonicity: existence, uniqueness, ergodicity, exponential mixing and  limit
theorems}

\author{{ \blfootnote{$^{*}$Corresponding author } Shuaishuai Lu$^{a}$ \footnote{ E-mail address : luss23@mails.jlu.edu.cn},~ Xue Yang$^{a}$}  \footnote{E-mail address : xueyang@jlu.edu.cn},~ Yong Li$^{a,b,*}$  \footnote{E-mail address : liyong@jlu.edu.cn}\\
	{$^{a}$College of Mathematics, Jilin University,} {Changchun 130012, P. R. China.}\\
	{$^{b}$School of Mathematics and Statistics, and Center for Mathematics and Interdisciplinary Sciences, Northeast Normal University,}
	{Changchun 130024, P. R. China.}
}

\begin{abstract}
This paper investigates the existence and uniqueness of solutions, as well as the ergodicity and exponential mixing to invariant measures, and limit theorems for a class of McKean-Vlasov SPDEs with locally weak monotonicity. In particular, for a class of weak monotonicity conditions, including H$\ddot{\text{o}}$lder continuity, we rigorously establish the existence and uniqueness of weak solutions to  McKean-Vlasov SPDEs by employing the Galerkin projection technique and the generalized coupling approach. Additionally, we explore the properties of the solutions, including time homogeneity, the Markov  and the Feller property. Building upon these properties, we examine the exponential ergodicity and mixing of invariant measures under Lyapunov conditions. Finally, within the framework of coefficients meeting the criteria of locally weak monotonicity and Lyapunov conditions, alongside the uniform mixing property of invariant measures, we establish the strong law of large numbers and the central limit theorem for the solution and obtain estimates of corresponding convergence rates.
~\\
~\\
 \textbf{keywords}: McKean-Vlasov SPDEs; Locally weak monotonicity; H$\ddot{\text{o}}$lder continuous  coefficients; Lyapunov condition; Exponential  ergodicity and mixing.
\end{abstract}
\end{frontmatter}
\section{\textup{Introduction}}
The McKean-Vlasov stochastic differential equations (SDEs), also called  mean-field SDEs, were initially elucidated by Kac \cite{ref6} within the framework of the Boltzmann equation addressing particle density in sparse monatomic gases and subsequently  was initiated by McKean \cite{ref7}. These equations find extensive utility across a spectrum of disciplines, notably in physics, statistical mechanics, quantum mechanics, and quantum chemistry.  For instance, as seen in the literature  \cite{ref8,ref9}, the utilization of McKean-Vlasov SDEs is demonstrated to characterize the asymptotic behavior of $N$-interacting particle systems with weak interactions manifesting chaotic propagation properties.  Larsy and Lions \cite{ref10,ref12}, as well as Huang et al. \cite{ref13,ref15}, respectively introduced mean-field strategies to investigate the deterministic behavior of large populations and stochastic differential strategies. Additionally, Michael \cite{ref16} analyzed the propagation of chaos results for the mean-field limit of a model for a trimolecular chemical reaction called the "Brusselator". Given the widespread applications of mean-field equations, we encourage interested readers to refer to \cite{ref17,ref19,ref20,ref21,ref22,ref23} and the references therein. This broad applicability stems from the recognition that the  evolution of many stochastic systems is contingent not solely upon the state of the solution but also upon its macroscopic distribution.

In the present paper, we will consider the following McKean-Vlasov stochastic (partial) differential equations:
\begin{align*}
		\begin{cases}
 \text{d}u(t)=(A(u(t), \mathcal{L}_{u(t)})+f(u(t),\mathcal{L}_{u(t)}))\text{d}t+g(u(t),\mathcal{L}_{u(t)})\text{d}W(t), \\
u(0) =x\in U_{1},
\end{cases}
	\end{align*}
 where $W(t)$ is a  cylindrical Wiener process  on a separable Hilbert space $(U_{2},\langle \cdot , \cdot \rangle_{ U_{2} })$ with respect to a complete
filtered probability space $\left ( \Omega ,\mathscr{F},\mathbb{P} \right )$, and $\mathcal{L}_{u(t)}$ denotes the law of $u(t)$ taking values in $\mathcal{P}(U_{1})$. When the system coefficients are Lipschitz continuous and satisfy linear growth conditions, there have been several relevant studies on the aforementioned  McKean-Vlasov stochastic systems. Wang \cite{ref25} established the  existence and uniqueness of solutions for McKean-Vlasov monotone SDEs, and   delved into pertinent topics such as exponential ergodicity and Harnack-type inequalities. Buckdahn et al. \cite{ref31}  demonstrated the correlation between functionals in the form of $\mathbb{E} f(t,X(t),\mathcal{L}_{x(t)})$ and the associated second-order PDEs and obtained the unique classical solution of a nonlocal mean-field  PDE.  In addition, under locally Lipschitz coefficients, Liu and  Ma \cite{ref3} derived conditions ensuring the global existence of solutions, subject to Lyapunov conditions. Furthermore, they  demonstrated the exponential convergence of invariant measures. For the McKean–Vlasov SPDEs with locally monotone conditions, Hong et al. \cite{ref111} investigated the strong and weak well-posedness. Barbu and R$\ddot{\text{o}} $ckner \cite{ref29} employed nonlinear Fokker-Planck equations to explore  the existence of weak solutions to McKean-Vlasov SDEs. For further results on this subject, we refer to \cite{ref32,ref35,ref36,ref57,ref58,ref65,ref27}, and others.

While the Lipschitz  continuity condition holds significant importance in the examination of dynamical properties within stochastic systems, it is commonly observed that the Lipschitz continuity condition is frequently violated by coefficients in numerous significant stochastic models. For instance, the Cox-Ingersoll-Ross model and the diffusion coefficients in the Ferrer branch diffusion  are merely H$\ddot{\text{o}}$lder continuous, rather than Lipschitz continuous. Consequently, for many such models, the (local) Lipschitz condition imposes a considerable constraint. As a result, stochastic models featuring non-Lipschitz coefficients have garnered growing interest in recent years. For example, Fang, Zhang \cite{ref37} and Wang et al. \cite{ref38} investigated SDEs and stochastic functional differential equations (SFDEs) featuring non-Lipschitz coefficients, specifically examining the existence and uniqueness of strong solutions.  In Reference \cite{ref50}, Kulik and Scheutzow utilized the method of generalized coupling to establish weak uniqueness for the weak solution of a specific class of SFDEs with H$\ddot{\text{o}}$lder continuous coefficients. Han \cite{ref999} established the well-posedness and small mass limit for the stochastic wave equation with H$\ddot{\text{o}}$lder noise coefficient. Furthermore, R$\ddot{\text{o}} $ckner et al. \cite{ref0000000}  investigated  the averaging principle of semilinear slow-fast stochastic partial differential equations with additive noise, employing the Poisson equation in Hilbert space under the assumption of H$\ddot{\text{o}}$lder continuity for the fast variables. However, the current results do not extend to stochastic partial differential systems in cases where the diffusion and drift coefficients only satisfy H$\ddot{\text{o}}$lder continuity. Hence, this paper will initially investigate the existence and uniqueness of (strong) solutions within the context of McKean-Vlasov SPDEs with locally weak monotonicity coefficients. And our framework will be appropriate for stochastic partial differential models with H$\ddot{\text{o}}$lder continuous coefficients.

As is known, the uniqueness plays a crucial role in determining the properties of (strong) solutions in stochastic differential systems. A classical finding reveals that under Lipschitz conditions, the influence of the initial value on the solution can be examined using Gronwall's lemma, which leads to the establishment of pathwise uniqueness for the weak solution. However, when the regularity of the coefficients falls below the Lipschitz condition, employing Gronwall's lemma becomes unfeasible. Consequently, one of the aims  in this paper is to address the aforementioned issues. Focusing on McKean-Vlasov S(P)DEs, we initially establish the existence and pathwise uniqueness of the weak solution for a  class of finite-dimensional systems, where the coefficients solely adhere to the conditions of  weak monotonicity and linear growth. We then proceed to establish the existence and uniqueness of the strong solution by leveraging the Yamada-Watanabe criterion. For detailed elaboration, refer to Theorem 3.1. In the context of infinite-dimensional McKean-Vlasov systems with coefficients satisfying the weak monotonicity, we employ the Galerkin projection technique. This method, combined with insights from finite-dimensional analysis, enables the proof of existence and uniqueness of variational solutions, as elaborated in Theorem 3.2. Under the weak monotonicity condition, where the diffusion and drift coefficients only satisfy H$\ddot{\text{o}}$lder continuity(indicating even lower regularity), the system's solution may exhibit localized oscillatory behavior, potentially leading to unmeasurable branching phenomena in the path space. Consequently, using traditional methods such as the contraction mapping principle or the martingale problem approach, the pathwise uniqueness of the weak solution is typically unattainable. As a result, the existence and uniqueness of a strong solution cannot be deduced from Yamada-Watanabe Theorem. Motivated by the generalized coupling method introduced in \cite{ref50}, we  investigate McKean-Vlasov SPDEs and establish the weak uniqueness of the solution under H$\ddot{\text{o}}$lder continuity, meaning that the weak solution is unique in law, as elaborated in Theorem 3.5.
	
Establishing limit theorems for Markov processes, particularly the  strong law of large numbers (SLLN) and the central Limit Theorem (CLT), pivotal in elucidating the long-term dynamics of stochastic processes, constitutes a central pursuit in probability theory, see e.g. \cite{ref51,ref52,ref53,ref55,ref56,ref59}. Hence, another primary objective of this paper is to extend these theorems to McKean-Vlasov SPDEs under weak monotonicity conditions. This endeavor necessitates an exploration of the exponential ergodicity. Previous investigations into the ergodic property of mean-field systems have primarily focused on scenarios characterized by monotonic or Lipschitz conditions, as evidenced by references \cite{ref60,ref61,ref62}, among others. Nonetheless, the exploration of lower regularity conditions, notably the H$\ddot{\text{o}}$lder continuous condition, remains relatively rare in the existing literature. This study endeavors to refine the findings established in reference \cite{ref3} and  make extensions to McKean-Vlasov SPDEs under weak monotonicity conditions. Drawing inspiration from references \cite{ref3} and \cite{ref63}, we employ Lyapunov functions and the It$\hat{\text{o}} $ formula to unveil the exponential convergence of invariant measures within the framework of the Wasserstein quasi-metric. Furthermore, by reinforcing the condition to encompass integrable Lyapunov functions, we ascertain the exponential mixing property of invariant measures in the Wasserstein metric. Further elaboration on these advancements is provided in Theorem 4.2. This outcome is simple and interesting in its own right, as it imposes constraints on Lyapunov functions rather than solely on system coefficients, thereby broadening the method's applicability.

As previously indicated, in the case of McKean-Vlasov SPDEs under weak monotonicity conditions, an invariant measure $\mu ^{*}$ exhibiting uniform mixing properties is derived. Subsequently, leveraging both the exponential mixing property  and the Markov property, we aim to demonstrate the applicability of the limit theorems delineated in \cite{ref5}  to the present system. Specifically, we establish the SLLN, the CLT, and their associated convergence rates, as detailed in Theorems 5.1 and 5.6:
 \begin{enumerate}[(1)]
		\item SLLN:
$$\frac{1}{t}\int_{0}^{t}\Phi (u(s;0,x))\text{d}s \to \int _{U_{1} } \Phi (x ) \mu ^{*}(\text{d}x ) \quad \text{as}\quad t\to \infty ,\quad \mathbb{P} -a.s.,$$
where $\Phi$ is the observation function;
        \item  CLT:
        $$\frac{1}{\sqrt{t} }\int_{0}^{t}[\Phi(u(s;0,x ))-\int _{\mathcal{H} } \Phi(x ) \mu ^{*}(\text{d}x )]\text{d}s \overset{W}{\to}\xi  ,$$
        where $\overset{W}{\to}$ means weak convergence and $\xi  $ is a normal random variable.
	\end{enumerate}

The paper is structured as follows. The subsequent section delineates key definitions, notation, and lemmas. In Section 3, we first rigorously establish the existence and  uniqueness of a strong solution for a specific class of finite-dimensional systems under weak monotonicity conditions. Subsequently, leveraging the sophisticated Galerkin projection technique and leveraging the insights garnered from the finite-dimensional context, we proceed to demonstrate the existence and uniqueness of solutions for the corresponding infinite-dimensional system. Subsequently, for a class of weak monotonicity conditions, including H$\ddot{\text{o}}$lder continuity, we rigorously establish the existence and weak uniqueness of weak solutions to  McKean-Vlasov SPDEs by employing the Galerkin projection technique and the generalized coupling approach. Section 4 elucidates  properties of the solution including time homogeneity, the Markov and the Feller property. Based on these properties, this section also establishes the ergodicity and exponential mixing of invariant measure under Lyapunov conditions. Following this, in Section 5, we focus on establishing the SLLN, the CLT, and estimating the corresponding convergence rates for McKean-Vlasov SPDEs.
	\section{\textup{Preliminaries}}
	The  $\left ( \Omega ,\mathscr{F},\mathbb{P} \right )$ be a certain complete probability space with a filtration $\lbrace\mathscr{F}_{t}\rbrace_{t\ge0}$ satisfying the usual condition. If $K$ is a matrix or a vector, $K^{'}$ is its transpose. For a matrix $K$, the norm is expressed as $\left \| K\right \|= \sqrt{\text{trace}(KK^{'})}$ and denote by $\left \langle \cdot , \cdot  \right \rangle $ the inner product of $ \mathbb{R}^{n} $. $(U_{i},\left \| \cdot  \right \| _{U_{i}})(i=1,2)$ are  separable Hilbert spaces with  inner product $ \left \langle \cdot, \cdot \right \rangle _{U_{i}}$ and   $(B,\left \| \cdot  \right \| _{B})$ is a reflexive Banach space such that $B\subset U_{1}$. $U_{i}^{*}$, $B^{*}$ are the dual spaces of $U_{i}$, $B$, and
\begin{align*}
		B\subset  U_{1}\subset B^{*},
	\end{align*}
where the embedding $B\subset U_{1}$ is continuous and dense,  hence
we have $U_{1}^{*} \subset B ^{*}$ continuously and densely. Let $_{B ^{\ast}} \langle \cdot , \cdot \rangle _{B} $  denote the dualization between $B^{\ast}$ and $B$, which shows that for all $x \in U_{1}$, $y \in B$,
\begin{align*}
		_{B ^{\ast}} \langle x, y\rangle _{B} = \langle x, y\rangle_{ U_{1}},
	\end{align*}
 and  ($B,U_{1}, B ^{\ast}$) is called a Gelfand triple.  For any $q \ge 1$  and  the Banach space $(Y, \left \|   \cdot  \right \|_{Y}   )$,  $\mathfrak{L} ^{q}(\Omega, Y) $ denotes the Banach space of all $Y$-value random variables:
\begin{eqnarray*}
		\mathfrak{L} ^{q}(\Omega, Y)=\left \{ y:\Omega \to Y :\mathbb{E}\left \| y \right \| ^{q}=\int _{\Omega } \left \| y \right \| ^{p}_{Y}\text{d}\mathbb{P}< \infty   \right \}.
	\end{eqnarray*}

Let  $\mathcal{B}(U_{1})$ denote the $\sigma  $-algebra generated by space $U_{1}$ and $\mathcal{P}(U_{1}) $ be the family of all probability measures  defined on $ \mathcal{B}(U_{1})$. We introduce the following Banach space $\mathcal{P} _{2}(U_{1})$ of signed measures $\mu$ on  $\mathcal{B}(U_{1})$ satisfying
\begin{align*}
		\left \| \mu  \right \| _{U_{1}}^{2}=\int _{U_{1}}(1+\left \| x  \right \| _{U_{1}  })^{2}\left | \mu  \right | (\text{d}x )<\infty,
	\end{align*}
where $\left|\mu\right| = \mu^{+} + \mu^{-}$ and $\mu = \mu^{+} - \mu^{-}$ is the Jordan decomposition of $\mu$. Let $\mathcal{P}^{*}(U_{1})=\mathcal{P}(U_{1})\cap\mathcal{P} _{2}(U_{1})  $ be the family of all probability measures on $(U_{1} , \mathcal{B}(U_{1} ))$ with the following metric:
 \begin{align*}
		d_{U_{1}}( \mu ,\nu ) :=\sup\{\left | \int \Phi  \text{d}\mu- \int \Phi  \text{d}\nu \right |:\left \| \Phi   \right \|_{BL}\le 1  \},
	\end{align*}
 where $\left \| \Phi   \right \|_{BL}:=\left \| \Phi   \right \| _{\infty }+Lip(\Phi  )$ and $\left \|\Phi  \right \| _{\infty }=\sup_{x\in U_{1} }\frac{|\Phi  (x)|}{(1+\left \|x  \right \|_{U_{1}} )^{2}} $, $Lip(\Phi )=\sup_{x\ne y}\frac{\left | \Phi (x)-\Phi  (y) \right | }{\left \|x-y  \right \|_{U_{1} }} $. Then it is not difficult to verify that the space  $(\mathcal{P}^{*}(U_{1}), d_{U_{1}})$ is a complete metric space; see
 \cite{ref1} for more details on $(\mathcal{P}^{*}(U_{1}), d_{U_{1}})$ and its properties.
 %where $\left \| F \right \|_{\mathcal{K}   L}:=\left \| F \right \| _{\infty }+Lip(F)$ and $\left \| F \right \| _{\infty }=\sup_{\varphi \in \mathcal{R} }\frac{|F(\varphi)|}{1+\left \| \varphi  \right \|^{r_{1}/2 }_{\mathcal{R} } } $, $Lip(F)=\sup_{\varphi_{1}\ne \varphi_{2}}\frac{\left | F(\varphi_{1})-F(\varphi_{2}) \right | }{d_{r_{1},r_{2}}(\varphi _{1},\varphi _{2}) } $ and $r_{1},r_{2}>0$,
% \begin{align*}
%		d_{r_{1},r_{2}}(\varphi _{1},\varphi _{2})=(1\wedge \left \| \varphi _{1}-\varphi _{2} \right \|^{r_{2}}_{\mathcal{R} } )\sqrt{1+\left \| \varphi _{1} \right \|^{r_{1}}_{\mathcal{R} }+\left \| \varphi _{2} \right \|^{r_{1}}_{\mathcal{R} } }, \quad \varphi _{1},\varphi _{2}\in \mathcal{R}.
%	\end{align*}
%Obviously, $(\varphi _{1},\varphi _{2})\mapsto d_{r_{1},r_{2}}(\varphi _{1},\varphi _{2})$ is a quasi-distance, i.e., it is symmetric, lower semi-continuous

Consider the following McKean-Vlasov stochastic partial differential equations(SPDEs):
\begin{align}\label{r1}
		\begin{cases}
 \text{d}u(t)=(A(u(t), \mathcal{L}_{u(t)})+f(u(t),\mathcal{L}_{u(t)}))\text{d}t+g(u(t),\mathcal{L}_{u(t)})\text{d}W(t), \\
u(0) =x\in U_{1},
\end{cases}
	\end{align}
 where $W(t)$ is a  cylindrical Wiener processes  on a separable Hilbert space $(U_{2},\langle \cdot , \cdot \rangle_{ U_{2} })$ with respect to a complete
filtered probability space $\left ( \Omega ,\mathscr{F},\mathbb{P} \right )$, and $\mathcal{L}_{u(t)}$ denotes the law of $u(t)$ taking values in $\mathcal{P}(U_{1})$. These are the measurable maps: $$A:B\times\mathcal{P}^{*}(U_{1})\to B^{*},\quad f:U_{1}\times\mathcal{P}^{*}(U_{1})  \to U_{1},\quad g:U_{1}\times\mathcal{P}^{*}(U_{1})  \to \mathscr{L}(U_{2},U_{1}),$$
where $\mathscr{L}(U_{2},U_{1})$ is  the space of all Hilbert-Schmidt operators from $U_{2}$ into $U_{1}$. The definition of variational solution to system \eqref{r1} is given as follows.
~\\
\\\textbf{Definition 2.1.}\cite{ref65} We call a continuous $U_{1}$-valued $\lbrace\mathscr{F}_{t}\rbrace_{t\ge0}$-adapted process $\{u(t)\}_{t\in[0,T]}$ is a solution of the system \eqref{r1}, if for its $\text{d}t\times \mathbb{P} $-equivalent class $\{\hat{u}(t) \}_{t\in[0,T]}$ satisfying $\hat{u}(t)\in \mathfrak{L} ^{p}([0,T]\times \Omega ,B)\times \mathfrak{L} ^{2}([0,T]\times \Omega ,U_{1})$
and $\mathbb{P}$-a.s.,
 \begin{align*}
		\text{d}u(t)=x+\int_{o}^{t} (A(\bar{u} (s), \mathcal{L}_{\bar{u} (s)})+f(\bar{u} (s),\mathcal{L}_{\bar{u} (s)}))\text{d}s+\int_{0}^{t} g(\bar{u} (s),\mathcal{L}_{\bar{u} (s)})\text{d}W(s),
	\end{align*}
where $\bar{u}$ is a $B$-valued progressively measurable $\text{d}t\times \mathbb{P} $ of $\hat{u}$.
~\\

The transition probability of the Markov process defined on $U_{1}$ is a function $p:\Delta \times U_{1} \times \mathcal{B}(U_{1} )\to\mathbb{R} ^{+} $, where $\Delta=\{(t,s):t\ge s, t,s \in \mathbb{R} \}$ for $u(t) \in U_{1} $ and $\Gamma\in \mathcal{B}(U_{1} )$ with the following properties:
\begin{enumerate}[1)]
		\item $p(t,s,x,\Gamma)=\mathbb{P}(\omega  :u(t;s,x)\in \Gamma|u_{s})$;
        \item $p(t,s,\cdot,\Gamma)$ is $\mathcal{B}(U_{1} )$-measurable for every $t\ge s$ and $\Gamma\in \mathcal{B}(U_{1} )$;
		\item $p(t,s,x ,\cdot)$ is a probability measure on $\mathcal{B}(U_{1} )$ for every $t\ge s$ and $x \in U_{1} $;
        \item The Chapman-Kolmogorov equation:
        \begin{align*}
		p(t,s,x,\Gamma)=\underset{U_{1} }{\int}p(t,\tau ,y,\Gamma)p(\tau,s ,x,\text{d}y )
	\end{align*}
holds for any $s \le \tau \le t  $, $x\in U_{1}$ and $\Gamma\in \mathcal{B}(U_{1} )$.
	\end{enumerate}
We further define a map $\hat{p} (t,s): \mathcal{P}(U_{1} )\to  \mathcal{P}(U_{1} )$ for any $\mu \in \mathcal{P}(U_{1} )$ and $\Gamma\in \mathcal{B}(U_{1} )$ by
\begin{align}\label{p23}
		\hat{p} (t,s)\mu(\Gamma)=\underset{U_{1} }{\int} p(t,s,y,\Gamma)\mu (\text{d}y).
	\end{align}
The Markov process $u(t)$ is said to be time-homogeneous with the transition function $p(t,s,x,\Gamma)$, if
\begin{align*}
		p(t ,s ,x,\Gamma)=p(t+\theta ,s+\theta ,x,\Gamma),
	\end{align*}
for all $t\ge s\ge0$, $\theta\ge0$, $x\in U_{1}$ and $ \Gamma \in \mathcal{B} (U_{1})$.

\section{\textup{Existence and uniqueness of solutions for McKean-Vlasov SPDEs under locally weak monotonicity conditions}}
To study the existence and uniqueness of  solutions for  \eqref{r1} under locally weak monotonicity conditions, we assume that the initial value $x\in U_{1}$ is independent of $W(t)$. We first assume that the coefficients in \eqref{r1} satisfy the following  hypotheses:
 ~\\
 \\\textbf{(H1)} (Continuity) For all $x,y \in B $ and $\mu \in\mathcal{P}^{*}(U_{1})$, the map
\begin{align*}
		B\times \mathcal{P}^{*}(U_{1})\ni (x,\mu )\to_{B ^{\ast}}\langle A(x,\mu),y\rangle _{B}
	\end{align*}
is continuous.
\\\textbf{(H2)} (Coercivity)  There exist constant  $ \lambda_{1}\in \mathbb{R} $, $\lambda _{2}>0 $, $p\ge2$ and $M>0$ such that for all $x \in B $, $\mu \in \mathcal{P}^{*}(U_{1})$
\begin{align*}
		_{B ^{\ast}} \langle A(x,\mu), x\rangle _{B} \le \lambda_{1}(\left \| x \right \| _{U_{1}}^{2}+\left \| \mu  \right \| _{U_{1}}^{2})-\lambda _{2}\left \| x \right \|_{B}^{p}+M.
	\end{align*}
 \\\textbf{(H3)} (Growth) For $A$, there exists constant $\lambda_{1}$ for all $x\in B$, $\mu \in \mathcal{P}^{*}(U_{1})$ such that
 \begin{align*}
		\left \| A(x,\mu)\right \|_{B^{*}}^{\frac{p}{p-1}}   \le \lambda  _{2}(\left \| x\right \|_{B}^{p}+\left \| \mu  \right \| _{U_{1}}^{2})+ M,
	\end{align*}
and for the continuous functions $f$, $g$, there exist constant $\lambda _{1}$ for all $x\in U_{1}  $, $\mu \in \mathcal{P}^{*}(U_{1}  )$ such that
\begin{align*}
		\left \| f(x,\mu) \right \|_{U_{1}} \vee \left \| g(x,\mu ) \right \| _{\mathscr{L}(U_{2},U_{1})}  \le \lambda _{2}(\left \| x\right \|_{U_{1}  }+\left \| \mu  \right \| _{U_{1} })+M.
	\end{align*}
\\\textbf{(H4)} (Locally weak monotonicity) Let $\mathcal{G}:[0,1)\to \mathbb{R} _{+} $ is an increasing, concave and continuous function satisfying $$\mathcal{G}(0)=0, \int_{0^{+}}\frac{\text{d}r}{\mathcal{G}(r)}=+\infty.$$ The map $A$ satisfies, for all $x,y\in B$, $\mu,\nu \in \mathcal{P}^{*}(U_{1})$,
\begin{align*}
		2_{B ^{\ast}} \langle A(x,\mu)-A(y,\nu), x-y\rangle _{B} \le \lambda_{1}[\mathcal{G}(\left \| x-y \right \| _{U_{1}}^{2})+d_{U_{1}}^{2} (\mu ,\nu)) ,
	\end{align*}
and   the functions $f$, $g$ satisfy,  for all $x,y\in U_{1}  $ and  $\mu, \nu\in \mathcal{P}^{*}(U_{1}  )$ with $\left \| x-y\right \| _{U_{1}}\vee d_{U_{1}} (\mu ,\nu)\le 1$,
\begin{align*}
		\left \langle f(x,\mu)-f(y,\nu),x-y \right \rangle _{U_{1}} \le\lambda_{1}[\mathcal{G}(\left \| x-y \right \| _{U_{1}}^{2})+d_{U_{1}}^{2} (\mu ,\nu)] ,
	\end{align*}
\begin{align*}
	 \left \| g(x,\mu )-g(y,\nu ) \right \|^{2} _{\mathscr{L}(U_{2},U_{1})} \le\lambda_{2}[\mathcal{G}(\left \| x-y \right \| _{U_{1}}^{2})+d^{2}_{U_{1}} (\mu ,\nu)].
	\end{align*}

In this paper, we write $L_{M}$ to mean some positive constants which depend on $M$.  If it does not cause confusion, we always assume that the constants $\lambda_{1}\in \mathbb{R} $, $\lambda_{2},M>0$ and these constants may change from line to line.  Now we discuss the existence and uniqueness of  solutions for  \eqref{r1} whose coefficients are assumed to be  Non-Lipschitz continuous.
~\\
\\\textbf{Theorem 3.1.} \emph{Consider \eqref{r1}. Suppose that the assumptions \textbf{(H1)}$-$\textbf{(H4)} hold. Assume also that $\lim_{r \to 0^{+}} \frac{\mathcal{G}(r) }{\sqrt{r} } =0.$, then for any initial values $x \in U_{1} $, system \eqref{r1} has a unique solution $u(t)$  in the sense of
Definition 2.1.}
~\\

Analyzing system \eqref{r1}, we will use the technique of Galerkin type approximation to get the existence and uniqueness of solutions. For ease of description, we first study the following McKean-Vlasov SDEs under locally weak monotonicity conditions:
\begin{align}\label{r3} \text{d}x(t)=F(x(t),\mathcal{L}_{x(t)})\text{d}t+G(x(t),\mathcal{L}_{x(t)})\text{d}B(t),
	\end{align}
with the initial data $x(0)=x_{0}\in \mathbb{R} ^{n}$, where $B(t)$ is an $m$-dimensional Wiener process and $F:\mathbb{R} ^{n}\times\mathcal{P}^{*}(\mathbb{R} ^{n})  \to \mathbb{R}^{n}$, $G:\mathbb{R} ^{n}\times\mathcal{P}^{*}(\mathbb{R} ^{n}) \to \mathbb{R}^{n\times m}$ are two continuous maps. To ensure the existence and uniqueness of solutions for \eqref{r3}, we make the following assumptions:
\\\textbf{(h1)} The functions $F$ and $G$ are continuous in $(x, \mu )$ and satisfy, for all $x \in \mathbb{R} ^{n}$, $\mu \in \mathcal{P}^{*}(\mathbb{R} ^{n})$
 \begin{eqnarray*}
		\left \langle F(x ,\mu ) ,x \right \rangle   \vee \left \| G(x, \mu  ) \right \|^{2} \le \lambda _{1}(\left \| x\right \|^{2}+\left \| \mu  \right \| _{\mathbb{R} ^{n}}^{2})+M.
	\end{eqnarray*}
\textbf{(h2)}  The functions $F,G$ satisfy, for all $x , y \in \mathbb{R} ^{n} $ and $\mu, \nu\in \mathcal{P}^{*}(\mathbb{R} ^{n} )$ with $\left \| x-y \right \| \le 1$ and $d_{\mathbb{R} ^{n}} (\mu,\nu)\le 1$,
 \begin{align*}
		2\left \langle F(x,\mu)-F(y,\nu),x-y\right \rangle \vee \left \| G(x,\mu )-G(y,\nu ) \right \|^{2} \le\lambda_{1}[\mathcal{G}(\left \| x-y \right \| _{U_{1}}^{2})+d^{2}_{\mathbb{R} ^{n}} (\mu,\nu) ].
	\end{align*}
~\\
\textbf{Theorem 3.2.} \emph{Consider \eqref{r3}. Suppose that the assumptions \textbf{(h1)}$-$\textbf{(h2)} hold. Then the following statement holds: for any $x_{0} \in \mathbb{R} ^{n} $, there exist a  unique strong solution $x(t)$  to \eqref{r3} with $x(0)=x_{0}$.
~\\
%\begin{enumerate}[1)]
%        \item  For any $\varphi \in \mathcal{R} $, there exist a  unique strong solution $u(t)_{t\in[0,T]}$ and segment process $u(t)$ to \eqref{r3} with $u_{0}=\varphi$;
%		\item  The segment process $u(t)$ is a time homogeneous Markov process in $\mathcal{R}$ and  has the Feller property.
%	\end{enumerate}
\\\textbf{proof:}} Based on the above analysis, we divide the proof of Theorem 3.2 into the following two steps:
\\\textbf{step 1:}
\textit{Existence of the weak solutions.} In the later proofs we will frequently consider Lipschitz approximations of  continuous mappings. By Proposition 3.3 of \cite{ref999},  we can find  sequences $\{F^{n}\},\{G^{n}\}$ satisfy the following properties:
\begin{enumerate}[(\textbf{a})]
		\item $F^{n}\to F$, $G^{n}\to G$ as $n\to \infty $  uniformly on each compact subset of $\mathbb{R} ^{n}\times\mathcal{P}^{*}(\mathbb{R} ^{n})$;
\end{enumerate}
\begin{enumerate}[(\textbf{b})]
        \item $F^{n},G^{n}$ satisfy the conditions \textbf{(h1)} and \textbf{(h2)}, that is,  the coefficients are independent of $n$.
	\end{enumerate}
\begin{enumerate}[(\textbf{c})]
        \item The functions $F^{n}$ and $G^{n}$ are Lipschitz continuous on each bounded subset of $\mathbb{R} ^{n}\times\mathcal{P}^{*}(\mathbb{R} ^{n})$.
	\end{enumerate}
 Fix $n\ge 1$ arbitrarily and we further consider the following equation:
\begin{align}\label{r5} \text{d}x^{n}(t)=F^{n}(x^{n}(t),\mathcal{L}_{x^{n}(t)})\text{d}t+G^{n}(x^{n}(t),\mathcal{L}_{x^{n}(t)})\text{d}B(t).
	\end{align}
From the above analysis, \eqref{r5}  has a unique strong solution $x^{n}(t)$ with the initial condition $x^{n}(0) =x_{0}$ by Theorem 2.1 of \cite{ref111}. Then for any fixed $q\ge2$, applying It$\hat{\text{o}} $ formula formula to $\|x^{n}(t)\|^{q}$, we have
\begin{align*}
\begin{split}
\|x^{n}(t)\|^{q}&=\left \| x_{0} \right \| ^{q}+\int_{0}^{t} [q\left \| x^{n}(s) \right \|^{q-2}\left \langle F^{n}(x^{n}(s),\mathcal{L}_{x^{n}(s)}),x^{n}(s)  \right \rangle \\&~~~+\frac{q(q-1)}{2} \left \| x^{n}(s) \right \|^{q-2} \left \| G^{n}(x^{n}(s),\mathcal{L}_{x^{n}(s)}) \right \|^{2}]\text{d}s
\\&~~~+\int_{0}^{t}q\left \| x^{n}(s) \right \|^{q-2} [x^{n}(s)]^{'} G^{n}(x^{n}(t),\mathcal{L}_{x^{n}(t)})\text{d}B(s).
\end{split}
\end{align*}
By  Young inequality and assumptions \textbf{(h1)} and \textbf{(h2)}, there exist some  constants $L_{q,\lambda _{1}}$ such that
\begin{align*}
\|x^{n}(t)\|^{q}&\le\left \| x_{0}  \right \|^{q}+L_{q,\lambda _{1}}\int_{0}^{t} (\left \| x^{n}(s) \right \|^{q}+\left \|\mathcal{L}_{x^{n}(s)}\right \| _{\mathbb{R} ^{n}}^{q}+M)\text{d}s
\\&~~~+L_{q,\lambda _{1}}\int_{0}^{t}(\left \| x^{n}(s) \right \| ^{q}+\left \|\mathcal{L}_{x^{n}(s)}\right \| _{\mathbb{R} ^{n}}^{q}+M)\text{d}B(s),
\end{align*}
where $L_{q,\lambda _{1}}$ depends only on $q,\lambda _{1}$.  By Cauchy’s inequality and Burkholder–Davis–Gundy inequality, we obtain
\begin{align*}
\begin{split}
\mathbb{E} \underset{z\in[0,t]}{\sup}  \left \| x^{n}(t) \right \| ^{2q}&\le 3 \left \| x_{0}  \right \|^{q}+L_{q,\lambda _{1}}\mathbb{E}[\int_{0}^{t} (\left \| x^{n}(s) \right \| ^{q}+\left \|\mathcal{L}_{x^{n}(s)}\right \| _{\mathbb{R} ^{n} }^{q}+M)\text{d}s]^{2}
\\&~~~+L_{q,\lambda _{1}}\mathbb{E}(\underset{z\in[0,t]}{\sup}\int_{0}^{t}(\left \| x^{n}(s) \right \| ^{q}+\left \|\mathcal{L}_{x^{n}(s)}\right \| _{\mathbb{R} ^{n}}^{q}+M)\text{d}B(s)) ^{2}
\\&\le 3 \left \| x_{0}  \right \|^{q}+L_{q,\lambda _{1}}\mathbb{E}\int_{0}^{t} (\left \| x^{n}(s) \right \| ^{2q}+\left \|\mathcal{L}_{x^{n}(s)}\right \| _{\mathbb{R} ^{n} }^{2q}+M^{2})\text{d}s
\\&\le3 \left \|x_{0} \right \| ^{q}+L_{q,\lambda _{1}}\int_{0}^{t} [\mathbb{E}\left \| x^{n}(s) \right \|^{2q}+L_{q}\mathbb{E}(1+\left \|x^{n}(s)\right \| ^{2q})+M^{2}]\text{d}s
\\&\le 3\left \| x_{0}  \right \|^{2q}+L_{q,\lambda _{1},M}[\int_{0}^{t} \mathbb{E}\underset{z\in[0,s]}{\sup}  \left \| x^{n}(z) \right \| ^{2q}\text{d}s+t].
\end{split}
\end{align*}
Then, applying the Gronwall inequality  yields
\begin{align}\label{r6}
\begin{split}
\mathbb{E} \underset{z\in[0,t]}{\sup}  \left \| x^{n}(t) \right \| ^{2q}\le L_{q,\lambda _{1},M}(\left \| x_{0}  \right \|^{2q}+t+e^{t})<\infty,
\end{split}
	\end{align}
for any $T>0$ and $t\in [0,T]$. We can follow from condition \textbf{(h1)} that $F^{n}$ and $G^{n}$ are bounded on every bounded subset of $\mathbb{R} ^{n}\times\mathcal{P}^{*}(\mathbb{R} ^{n})$, further by \textbf{(h1)} and \eqref{r6}, there exists a constant $L_{q,\lambda _{1},M,T}>0$ independent of $n$ such that
\begin{align*}
\left \|F^{n}(x^{n}(s),\mathcal{L}_{x^{n}(s)} )\right \| \vee \left \|  G^{n}(x^{n}(s),\mathcal{L}_{x^{n}(s)} )\right \|\le L_{q,\lambda _{1},M,T}.
\end{align*}
Hence for any $0\le z, t \le T < \infty  $, we have
\begin{align}\label{r7}
\begin{split}
&\underset{n\ge1}{\sup}\mathbb{E} \left \| x^{n}(t)-x^{n}(z) \right \| ^{2q}\\&\le L_{q,\lambda _{1}}\underset{n\ge1}{\sup}\mathbb{E} \left \| \int_{z}^{t} F^{n}(x^{n}(s),\mathcal{L}_{x^{n}(s)})\text{d}s \right \| ^{2q}
+L_{q,\lambda _{1}}\mathbb{E} \underset{n\ge1}{\sup}\left \| \int_{z}^{t} G^{n}(x^{n}(s),\mathcal{L}_{x^{n}(s)})\text{d}B(s) \right \| ^{2q}
\\&\le L_{q,\lambda _{1},M,T}\left \| t-z \right \| ^{q},
\end{split}
	\end{align}
which yields that the family of laws $\mathcal{L}_{x^{n}(t)}$ of $x^{n}(t)$  is weakly compact. Further, from property \textbf{(a)} we can obtain the weak limit point $x^{*}(t)$ of $x^{n}(t)(n\to \infty)$, being a weak solution to \eqref{r3}. This proof is again quite standard and see \textbf{Appendix \uppercase\expandafter{\romannumeral1}} for details.
~\\
\\\textbf{step 2:}
\textit{Pathwise uniqueness of the weak solutions.} In the following, we will present the pathwise uniqueness for \eqref{r3}. Suppose that two stochastic processes $x(t),y(t)$ satisfy the following form:
\begin{align*}
 x(t)=x_{0}+\int_{0}^{t} F(x(s),\mathcal{L}_{x(s)})\text{d}s+\int_{0}^{t} G(x(s),\mathcal{L}_{x(s)})\text{d}B(s),
	\end{align*}
and
\begin{align*}	
 y(t)=y_{0}+\int_{0}^{t} F(y(s),\mathcal{L}_{y(s)})\text{d}s+\int_{0}^{t} G(y(s),\mathcal{L}_{y(s)})\text{d}B(s).
	\end{align*}
Without loss of generality, let's assume $\left \| x_{0} -y_{0} \right \|<1$. Denote
\begin{align*}
		\tau _{N}=\underset{t\ge0}{\inf} \{\left \| x(t)\right \|\vee \left \|y(t)\right \|>N \},\quad \tau _{\gamma }=\underset{t\ge0}{\inf} \{\left \| x(t)-y(t)\right \|>\gamma  \},
	\end{align*}
where $N>\left \| x_{0}  \right \|\vee \left \| y_{0}  \right \|$ and $\gamma \in(\left \| x_{0} -y_{0}   \right \|,1 ]$. In order to better prove the pathwise uniqueness of the weak solution, we first assume that the following property holds:
\begin{align}\label{r12}
		 \lim_{\left \| x_{0}  -y_{0}   \right \| \to 0} \mathbb{E} (\underset{z\in[0,t]}{\sup}\left \| x(z\wedge \tau _{N }\wedge \tau _{\gamma } )-y(z\wedge \tau _{N }\wedge \tau _{\gamma } ) \right \| ^{2})=0.
	\end{align}
Then by \eqref{r6}, we obtain that $\lim_{N \to \infty} \tau _{N}=\infty $ a.s. Then by  Fatou's lemma,  if $\left \| x_{0} -y_{0}   \right \| =0$, we obtain  $$ \mathbb{E} (\underset{z\in[0,t]}{\sup}\left \| x(z\wedge \tau _{\gamma } )-y(z\wedge \tau _{\gamma } ) \right \| ^{2})=0.$$ In addition, by the definition of $\tau _{\gamma }$, we have $\left \| x(\tau _{\gamma })- y(\tau _{\gamma }) \right \|=\gamma  $, which implies
\begin{align*}
		 \mathbb{E} (\underset{z\in[0,t]}{\sup}\left \| x(z\wedge \tau _{\gamma } )-y(z\wedge \tau _{\gamma } ) \right \| ^{2})&=\mathbb{P}(t\ge\tau _{\gamma})\mathbb{E} (\underset{z\in[0,t]}{\sup}\left \| x(z\wedge \tau _{\gamma } )-y(z\wedge \tau _{\gamma } ) \right \| ^{2})
\\&~~~+\mathbb{P}(t<\tau _{\gamma})\mathbb{E} (\underset{z\in[0,t]}{\sup}\left \| x(z\wedge \tau _{\gamma } )-y(z\wedge \tau _{\gamma } ) \right \| ^{2})
\\&\ge\mathbb{P}(t\ge\tau _{\gamma}) \gamma^{2},
	\end{align*}
i.e., $\mathbb{P}(t\ge\tau _{\gamma})=0$ as $\left \| x_{0} -y_{0}  \right \| \to 0$. Hence, by  \textbf{(h1)} and \eqref{r6}, we have
\begin{align*}
		 \mathbb{E} (\underset{z\in[0,t]}{\sup}\left \| x(z)-y(z) \right \| ^{2})&=\mathbb{E} (\underset{z\in[0,t]}{\sup}\left \| x(z )-y(z ) \right \| ^{2}\chi_{\{t>\tau _{\gamma }\}})
+\mathbb{E} (\underset{z\in[0,t]}{\sup}\left \| x(z )-y(z ) \right \| ^{2}\chi_{\{t\le\tau _{\gamma }\}})
\\&\le \mathbb{E} (\underset{z\in[0,t]}{\sup}\left \| x(z\wedge \tau _{\gamma } )-y(z\wedge \tau _{\gamma } ) \right \| ^{2})\\&=0,
	\end{align*}
which implies that as $\left \| x_{0}-y_{0} \right \| \to 0$,
\begin{align*}
		\mathbb{P}(\underset{z\in[0,t]}{\sup}\left \| x(z)-y(z) \right \| ^{2}=0) =1,
	\end{align*}
i.e., the  pathwise uniqueness of weak solution holds, and then by Yamada–Watanabe principle, we obtain that there exists a unique global strong solution of \eqref{r3}.

From the above analysis, we know that the key condition for pathwise uniqueness  is that \eqref{r12} holds.
Then applying the It$\hat{\text{o}} $ formula to $ \left \| x(t\wedge \tau _{N }\wedge \tau _{\gamma } )-y(t\wedge \tau _{N }\wedge \tau _{\gamma } ) \right| ^{2} $, Young inequality, BDG inequality and \textbf{(h2)} yield
\begin{align}\label{r8}
		 & \mathbb{E} (\underset{z\in[0,t]}{\sup}\left \| x(z\wedge \tau _{N }\wedge \tau _{\gamma } )-y(z\wedge \tau _{N }\wedge \tau _{\gamma } ) \right \| ^{2}) \nonumber
\\&=\left \|  x_{0}-y_{0}   \right \| ^{2}+\mathbb{E}\int_{0}^{t\wedge \tau _{N }\wedge \tau _{\gamma }} [\left \| G(x(s),\mathcal{L}_{x(s)})-G(y(s),\mathcal{L}_{y(s)}) \right \| ^{2} \nonumber
  \\&~~~+2\left \langle F(x(s),\mathcal{L}_{x(s)})-F(y(s),\mathcal{L}_{y(s)}),x(s)-y(s) \right \rangle]\text{d}s      \nonumber
\\&~~~+2\mathbb{E} (\underset{z\in[0,t]}{\sup}\int_{0}^{z\wedge \tau _{N }\wedge \tau _{\gamma }}[x(s)-y(s)]'( G(x(s),\mathcal{L}_{x(s)})-G(y(s),\mathcal{L}_{y(s)}) )\text{d}B(s))
\\&\le\left \|  x_{0}-y_{0}  \right \| ^{2}+2\lambda _{1}\mathbb{E}\int_{0}^{t\wedge \tau _{N }\wedge \tau _{\gamma }} (\mathcal{G} (\left \| x(s)-y(s) \right \| ^{2})+d _{\mathbb{R} ^{n}}^{2} (\mathcal{L}_{x(s)},\mathcal{L}_{y(s)} ) )\text{d}s    \nonumber
\\&~~~+12\mathbb{E}[\int_{0}^{t\wedge \tau _{N }\wedge \tau _{\gamma }} \left \|(x(s)-y(s))'( G(x(s),\mathcal{L}_{x(s)})-G(y(s),\mathcal{L}_{y(s)}) )\right \| ^{2}\text{d}s]^{\frac{1}{2} }   \nonumber
\\&\le\left \|  x_{0}-y_{0}  \right \| ^{2}+2\lambda _{1}\mathbb{E}\int_{0}^{t\wedge \tau _{N }\wedge \tau _{\gamma }} (\mathcal{G} (\left \| x(s)-y(s) \right \|^{2 })+d _{\mathbb{R} ^{n}}^{2 } (\mathcal{L}_{x(s)},\mathcal{L}_{y(s)} ))\text{d}s    \nonumber
\\&~~~+\frac{1}{2} \mathbb{E} (\underset{z\in[0,t\wedge \tau _{N }\wedge \tau _{\gamma }]}{\sup}\left \| x(z)-y(z) \right \|^{2} )+72\lambda _{1}\mathbb{E}\int_{0}^{t\wedge \tau _{N }\wedge \tau _{\gamma }} [\mathcal{G} (\left \| x(s)-y(s) \right \| ^{2 })+d _{\mathbb{R} ^{n}}^{2 } (\mathcal{L}_{x(s)},\mathcal{L}_{y(s)} )]\text{d}s.  \nonumber
	\end{align}
By the definitions of $d_{\mathbb{R} ^{n}}(\cdot ,\cdot )$ ,
\begin{align}\label{r9}
		d _{\mathbb{R} ^{n}} (\mathcal{L}_{x(s)},\mathcal{L}_{y(s)} )&=\sup\{\left \| \int_{\mathbb{R} ^{n}} \Phi  \text{d}\mathcal{L}_{x(s)}- \int_{\mathbb{R} ^{n}} \Phi  \text{d}\mathcal{L}_{y(s)} \right \|:\left \| \Phi    \right \|_{BL}\le 1  \}   \nonumber
\\&=\sup\{\left \| \int_{\mathbb{R} ^{n}} \Phi (x ) \mathbb{P}(\omega :x(s)\in\text{d}x )- \int_{\mathbb{R} ^{n}} \Phi (x ) \mathbb{P}(\omega :y(s)\in\text{d}x  ) \right \|:\left \| \Phi    \right \|_{BL}\le 1  \}  \nonumber
\\&=\sup\{\left \| \mathbb{E} \Phi  (x(s))-\mathbb{E} \Phi  (y(s)) \right \|:\left \| \Phi    \right \|_{BL}\le 1  \}
\\&\le \sup\{\left \| Lip( \Phi   )\right \|\cdot \mathbb{E} \left \|x (s)-y(s) \right \|_{\mathbb{R} ^{n}}:\left \| \Phi    \right \|_{BL}\le 1  \} \nonumber
\\&\le\mathbb{E} \left \|x (s)-y(s) \right \|.\nonumber
	\end{align}
Substituting \eqref{r9}  into \eqref{r8}, by H$\ddot{\text{o}} $lder inequality and Jensen's Inequality we have
\begin{align}\label{r11}
		 & \mathbb{E} (\underset{z\in[0,t]}{\sup}\left \| x(z\wedge \tau _{N }\wedge \tau _{\gamma } )-y(z\wedge \tau _{N }\wedge \tau _{\gamma } ) \right \| ^{2}) \nonumber
\\&\le2\left \| x_{0}-y_{0} \right \| ^{2}+\lambda _{1}\mathbb{E}\int_{0}^{t\wedge \tau _{N }\wedge \tau _{\gamma }} \mathcal{G} (\underset{z\in[0,s]}{\sup}\left \| x(s)-y(s) \right \|^{2})+(\underset{z\in[0,s]}{\sup}\left \| x(z)-y(z ) \right \|^{2 }) \text{d}s
\\&\le 2\left \| x_{0}-y_{0} \right \| ^{2}+\lambda _{1}\int_{0}^{t} \mathcal{G} (\mathbb{E}\underset{z\in[0,s\wedge \tau _{N }\wedge \tau _{\gamma }]}{\sup}\left \| x(s)-y(s) \right \|^{2})+\mathbb{E}(\underset{z\in[0,s\wedge \tau _{N }\wedge \tau _{\gamma }]}{\sup}\left \| x(z)-y(z ) \right \|^{2 }) \text{d}s\nonumber
\\&:=\Lambda (t). \nonumber
	\end{align}
For any fix $\iota>0$, let $H(t)=\int_{\iota}^{t}(\mathcal{G}(s)+s)^{-1} \text{d}s  $, then
\begin{enumerate}[(\textbf{a})]
		\item $H(t)$ is a monotonically increasing function;
\end{enumerate}
\begin{enumerate}[(\textbf{b})]
        \item $\lim_{t \to 0^{+}} H(t)=-\infty$ ,
	\end{enumerate}
i.e., $H(t)$ satisfies $H(t)>-\infty$ for any $t > 0$. Hence
\begin{align}\label{r13}
		H(\mathbb{E} (\underset{z\in[0,t]}{\sup}\left \| x(z\wedge \tau _{N }\wedge \tau _{\gamma } )-y(z\wedge \tau _{N }\wedge \tau _{\gamma } ) \right \| ^{2}))\le H(\Lambda (t)),
	\end{align}
further
\begin{align}\label{r15}
		&H(\Lambda (t))\nonumber
\\&=H(\Lambda (0))+\int_{0}^{t}H'(\Lambda (s))\text{d}\Lambda (s) \nonumber
\\&=H(2\left \|  x_{0}-y_{0} \right \| ^{2})+\lambda _{1}\int_{0}^{t}\frac{\mathcal{G} (\mathbb{E}\underset{z\in[0,s\wedge \tau _{N }\wedge \tau _{\gamma }]}{\sup}\left \| x(s)-y(s) \right \|^{2})+\mathbb{E}(\underset{z\in[0,s\wedge \tau _{N }\wedge \tau _{\gamma }]}{\sup}\left \| x(z)-y(z ) \right \|^{2 }) }{ \mathcal{G}(\Lambda (s))+\Lambda (s) }\text{d}s
\\&\le H(2\left \|  x_{0}-y_{0} \right \| ^{2})+\lambda _{1} t.\nonumber
	\end{align}
By \eqref{r13} and \eqref{r15}, we obtain when $\left \|  x_{0}-y_{0} \right \|\to 0$, we have $$H(\mathbb{E} (\underset{z\in[0,t]}{\sup}\left \| x(z\wedge \tau _{N }\wedge \tau _{\gamma } )-y(z\wedge \tau _{N }\wedge \tau _{\gamma } ) \right \| ^{2}))\to -\infty	,$$ which implies
\begin{align*}
	\lim_{\left \| x_{0}  -y_{0}   \right \| \to 0} \mathbb{E} (\underset{z\in[0,t]}{\sup}\left \| x(z\wedge \tau _{N }\wedge \tau _{\gamma } )-y(z\wedge \tau _{N }\wedge \tau _{\gamma } ) \right \| ^{2})=0.
	\end{align*}
  This completes
the proof. \quad $\Box$
~\\
\\Now we give the proof of \textbf{Theorem 3.1.}
~\\
\\\textbf{Proof of Theorem 3.1.} We first apply the Galerkin projection technique to transform the system \eqref{r1} into a finite-dimensional system. Assuming an orthonormal basis $\left \{ \varepsilon _{1},\varepsilon _{2},\varepsilon _{3} ,...\right \} \subset B$ on $U_{1}$ and $\left \{ \iota  _{1},\iota  _{2},\iota  _{3} ,...\right \}$ on $U_{2}$ and taking the first $k$ orthonormal bases yield the following operators:
\begin{align*}
\begin{split}
		\mathscr{U} _{k}:B^{*}\to U_{1}^{k}:=\text{span}\{\varepsilon _{1},\varepsilon _{2},...,\varepsilon _{k} \},\quad \mathscr{W}  _{k}:U_{2}\to U_{2}^{k}:=\text{span}\{\iota _{1},\iota _{2},...,\iota _{k} \}.
\end{split}
	\end{align*}
Let
\begin{align*}
		u^{k}:=\mathscr{U}_{k}(u)=\sum_{i=1}^{k}~_{B ^{\ast}} \langle u, \varepsilon _{i}\rangle _{B}\varepsilon _{i},\quad W^{k}(t):=\mathscr{W}_{k}[W(t)]=\sum_{i=1}^{k}\left \langle W(t),\iota _{i} \right \rangle _{U_{2}}\iota _{i},
	\end{align*}
where $u\in B$ and $k\ge1$. We thus obtain the following finite-dimensional equations corresponding to system \eqref{r1}:
\begin{align}\label{r2}
		\begin{cases}
 \text{d}u^{k}(t)=\mathscr{U}_{k}[(A(u(t), \mathcal{L}_{u(t)})+f(u(t),\mathcal{L}_{u(t)}))]\text{d}t+\mathscr{U}_{k}[g(u(t),\mathcal{L}_{u(t)})]\text{d}W^{k}(t), \\
u^{k}(s) =x^{k}\in  U^{k}_{1},
\end{cases}
	\end{align}
where $u^{k}(t)=\mathscr{U}_{k}[u(t)]$ and $x^{k}:=\mathscr{U}_{k}(x)$. Following Theorem 3.2, by \textbf{(H1)},\textbf{(H3)} and \textbf{(H4)}, system \eqref{r2} has a unique strong solution $u^{k}(t)$. Next, we prove Theorem 3.2. Similarly, we prove it in the following two steps:
\\\textbf{step 1:}
\textit{Apriori estimates of the  solutions $u^{k}(t)$.} By It\^o's formula, \textbf{(H2)} and \textbf{(H3)}, we have
\begin{align}\label{p1}
		\left \| u^{k}(t) \right \|^{2}_{U_{1}} &=\left \| x ^{k} \right \|^{2}_{U_{1} } +\int_{0}^{t}[2_{B ^{\ast}} \langle \mathscr{U}_{k}(A(u^{k}(s), \mathcal{L}_{u^{k}(s)}), u^{k}(s)\rangle _{B}\nonumber
\\&~~~+2 \langle \mathscr{U}_{k}(f(u^{k}(s), \mathcal{L}_{u^{k}(s)}), u^{k}(s)\rangle _{U_{1}} +\left \|  \mathscr{U}_{k}[g(u^{k}(s),\mathcal{L}_{u^{k}(s)})]\mathscr{W}_{k}\right \|^{2}_{\mathscr{L}(U_{2},U_{1})} ]\text{d}s\nonumber
\\&~~~+2\int _{0}^{t}\langle u^{k}(s),\mathscr{U}_{k}[g(u^{k}(s),\mathcal{L}_{u^{k}(s)})]\text{d}W^{k}(s)\rangle_{U_{1}}
\\&\le\left \| x ^{k} \right \|^{2}_{U_{1} }+\int_{0}^{t}[-\lambda _{2}\left \| u^{k}(s) \right \|_{B}^{p}+\lambda_{1}(\left \| u^{k}(s) \right \| _{U_{1}}^{2}+\left \| \mathcal{L}_{u^{k}(s)} \right \| _{U_{1}}^{2}\nonumber
)+M]\text{d}s
\\&~~~+2\int _{0}^{t}\langle u^{k}(s),\mathscr{U}_{k}[g(u^{k}(s),\mathcal{L}_{u^{k}(t)})]\text{d}W^{k}(s)\rangle_{U_{1}}.\nonumber
	\end{align}
Denote
\begin{align*}
		\tau ^{k}_{N}=\underset{t\ge0}{\inf} \{\left \| u^{k}(t)\right \|_{U_{1}}>N \}.
	\end{align*}
Due to $\lambda _{2}>0$, according  to  the B-D-G inequality, Young's inequality, the definitions of $\left \| \mu  \right \| _{U_{1} }$, we get
\begin{align}\label{p2}
		&\mathbb{E} \underset{z\in[0,T\wedge \tau ^{k}_{N}]}{\sup}\left \| u^{k}(z) \right \|^{2}_{U_{1}} \nonumber
\\&\le\lambda_{1}\left \| x ^{k} \right \|^{2}_{U_{1} }+\lambda_{1}\mathbb{E}\int_{0}^{T\wedge \tau ^{k}_{N}}\underset{z\in[0,s]}{\sup}\left \| u^{k}(z) \right \| _{U_{1}}^{2}\text{d}s+MT \nonumber
\\&~~+6\mathbb{E}\int _{0}^{T\wedge \tau ^{k}_{N}} \left \| u^{k}(t) \right \| ^{2}_{U_{1}}\left \|\mathscr{U}_{k}[g(u^{k}(s),\mathcal{L}_{u^{k}(t)})]\right \| ^{2}_{\mathscr{L}(U_{2},U_{1})} \text{d}s
\\&\le\lambda_{1}\left \| x ^{k} \right \|^{2}_{U_{1} }+\lambda_{1}\mathbb{E}\int_{0}^{T\wedge \tau ^{k}_{N}}\underset{z\in[0,s]}{\sup}\left \| u^{k}(z) \right \| _{U_{1}}^{2}\text{d}s+MT \nonumber
\\&~~~+\frac{1}{2}\mathbb{E} \underset{z\in[0,T\wedge \tau ^{k}_{N}]}{\sup}\left \| u^{k}(t) \right \|^{2}_{U_{1}} \nonumber.
	\end{align}
By Gronwall's lemma, we have
\begin{align}\label{p3}
		\mathbb{E} \underset{z\in[0,T\wedge \tau ^{k}_{N}]}{\sup}\left \| u^{k}(z) \right \|^{2}_{U_{1}} &\le M_{T,\lambda_{1}}(1+\left \| x \right \|^{2}_{U_{1} }).
	\end{align}
Take expectations on both sides of \eqref{p1}.  Then by \eqref{p3} we obtain
\begin{align}\label{p5}
		\mathbb{E}\int _{0}^{T\wedge \tau ^{k}_{N}}\left \| u^{k}(s) \right \|_{B}^{p}\text{d}s\le M_{T,\lambda_{1}}(1+\left \| x \right \|^{2}_{U_{1} }).
	\end{align}
 Combining \textbf{H2}, \textbf{H3}, \eqref{p3} and \eqref{p5} yields
\begin{align}\label{p6}
		&\mathbb{E}\int _{0}^{T\wedge \tau ^{k}_{N}}[\left \| (A(u^{k}(s), \mathcal{L}_{u^{k}(s)}) \right \|_{B^{*}}^{\frac{p}{p-1}} +\left \| f(u^{k}(s), \mathcal{L}_{u^{k}(s)}) \right \|^{2}_{U_{1}}\nonumber
\\&~~~+\left \| g(u^{k}(s), \mathcal{L}_{u^{k}(s)} ) \right \| ^{2}_{\mathscr{L}(U_{2},U_{1})} ]\text{d}s
\\&\le M_{T,\lambda_{1}}(1+\left \| x \right \|^{2}_{U_{1} }).\nonumber
	\end{align}
Letting $N\to\infty $, we have $\tau ^{k}_{N}\to \infty $ by \eqref{r6}. We thus  obtain the following  priori estimates for the solution $u^{k}(s)$ of \eqref{r2}:
\begin{align}\label{p7}
		&\mathbb{E} \underset{z\in[0,T]}{\sup}\left \| u^{k}(z) \right \|^{2}_{U_{1}} +\mathbb{E}\int _{0}^{T}[\left \| (A(u^{k}(s), \mathcal{L}_{u^{k}(s)}) \right \|_{B^{*}}^{\frac{p}{p-1}}\nonumber
+\left \| f(u^{k}(s), \mathcal{L}_{u^{k}(s)}) \right \|^{2}_{U_{1}}
 \\&+\left \| g(u^{k}(s), \mathcal{L}_{u^{k}(s)} ) \right \| ^{2}_{\mathscr{L}(U_{2},U_{1})}+\left \| u^{k}(s) \right \|_{B}^{p} ]\text{d}s
\\&\le M_{T,\lambda_{1}}(1+\left \| x \right \|^{2}_{U_{1} }) \nonumber.
	\end{align}
\textbf{step 2:}
\textit{Existence and uniqueness of the solution to system \eqref{r1}.} According to the reflexivity of $\left \| \cdot\right \|_{B}^{p}$ and step 1, we may assume that there exist common subsequences $k_{L}$ such that when $L\to \infty  $:
\begin{enumerate}[\textbf{1})]
		\item $u^{k_{L}}(t)\to u(t)$ in $\mathfrak{L} ^{2}([0,T]\times \Omega ,U_{1})$ and  weakly in $\mathfrak{L} ^{p}([0,T]\times \Omega ,B)$;
\end{enumerate}
\begin{enumerate}[\textbf{2})]
        \item $A(u^{k_{L}}(\cdot ),\mathcal{L}_{u^{k}(\cdot )} )\to A^{*}$ weakly in $[\mathfrak{L} ^{p}([0,T]\times \Omega ,B)]^{*}$ and $f(u^{k_{L}}(\cdot ),\mathcal{L}_{u^{k_{L}}(\cdot )} )\to f^{*}$ weakly in $\mathfrak{L} ^{2}([0,T]\times \Omega ,U_{1})$;
	\end{enumerate}
\begin{enumerate}[\textbf{3})]
        \item $g(u^{k_{L}}(\cdot ),\mathcal{L}_{u^{k_{L}}(\cdot )} )\to g^{*}$ weakly in $\mathfrak{L} ^{2}([0,T]\times \Omega ,\mathscr{L}(U_{2},U_{1}))$ and hence
            $$\int_{0}^{t} g(u^{k_{L}}(s),\mathcal{L}_{u^{k_{L}}(s)})\text{d}W(s)\to \int_{0}^{t}g^{*}(s)\text{d}W(s)$$
            weakly* in $L^{\infty }([0,T], L^{2}(\Omega , U_{1}))$.
	\end{enumerate}
Note that $B$ is separable, and hence  for any $v \in B  $ and $t\in[0,T]$, we obtain
\begin{align*}
		&\mathbb{E} \int_{0}^{t}  [_{B ^{\ast}} \langle u(s), v \rangle _{B}]\text{d}s
\\&=\lim_{L \to \infty} \mathbb{E} \int_{0}^{t}  [_{B ^{\ast}} \langle u^{k_{L}}(s), v \rangle _{B}]\text{d}s
\\&=\lim_{L \to \infty}\mathbb{E}\int_{0}^{t}[_{B ^{\ast}} \langle x^{k_{L}} , v \rangle _{B}+\int_{0}^{s}  (_{B ^{\ast}} \langle A(u^{k_{L}}(z),\mathcal{L}_{u^{k_{L}}(z)} ), v \rangle _{B})\text{d}z
\\&~~~+\int_{0}^{s}\left \langle f(u^{k_{L}}(z),\mathcal{L}_{u^{k_{L}}(z)} ),v \right \rangle _{U_{1}}\text{d}z
+\int _{0}^{s}\langle v,g(u^{k_{L}}(z),\mathcal{L}_{u^{k_{L}}(z)})\text{d}W(z)\rangle_{U_{1}}]\text{d}s
\\&=\mathbb{E}\int_{0}^{t}[_{B ^{\ast}} \langle x , v \rangle _{B}+\int_{0}^{s}  (_{B ^{\ast}} \langle A^{*}(z), v \rangle _{B})\text{d}z
+\int_{0}^{s}\left \langle f^{*}(z),v \right \rangle _{U_{1}}\text{d}z
\\&~~~+\int _{0}^{s}\langle v,g^{*}(z)\text{d}W(z)\rangle_{U_{1}}]\text{d}s,
	\end{align*}
which implies for any $t\in [0,T]$,
\begin{align*}
u(t)=x+\int_{0}^{t}A^{*}(s)\text{d}s+\int_{0}^{t}f^{*}(s)\text{d}s+\int_{0}^{t}g^{*}(s)\text{d}W(s), \quad \text{d}t\times \mathbb{P}-\text{a.e.}
	\end{align*}
Thus, it suffices to prove that
\begin{align}\label{p8}
A^{*}=A(u(\cdot  ),\mathcal{L}_{u(\cdot  )}), \quad f^{*}=f(u(\cdot  ),\mathcal{L}_{u(\cdot  )}),\quad g^{*}=g(u(\cdot  ),\mathcal{L}_{u(\cdot  )}), \quad \text{d}t\times \mathbb{P}-\text{a.e.}
	\end{align}
In fact, give any $\varrho (t) \in \mathfrak{L} ^{p}([0,T]\times \Omega ,B)\cap \mathfrak{L} ^{2}([0,T]\times \Omega ,U_{1})$. Without loss of generality, we assume that $\mathbb{E}\left \| u^{k_{L}}(s)-\varrho(s) \right \| _{U_{1}}^{2}\le1$ for any $t\in [0,T]$. Let $\eta =
\left | \lambda _{1} \right | $, then applying the It\^o's formula, we obtain
\begin{align}\label{p9}
		 &\mathbb{E}e^{-\eta  t}\left \| u^{k_{L}}(t) \right \|_{U_{1}}^{2}  -\left \|x  \right \|_{U_{1}}^{2}   \nonumber
\\&=\mathbb{E} \int_{0}^{t}e^{-\eta  s}[2_{B ^{\ast}} \langle A(u^{k_{L}}(s),\mathcal{L}_{u^{k_{L}}(s)} ), u^{k_{L}}(s) \rangle _{B}+2\langle f(u^{k_{L}}(s),\mathcal{L}_{u^{k_{L}}(s)} ), u^{k_{L}}(s) \rangle _{U_{1}}  \nonumber
\\&~~~+\left \| g(u^{k_{L}}(s),\mathcal{L}_{u^{k_{L}}(s)} ) \right \|_{\mathscr{L}(U_{2},U_{1})}^{2}- \eta\left \| u^{k_{L}}(s) \right \|_{U_{1}}^{2}]\text{d}s  \nonumber
\\&\le \mathbb{E} \int_{0}^{t}e^{-\eta  s}[2_{B ^{\ast}} \langle A(u^{k_{L}}(s),\mathcal{L}_{u^{k_{L}}(s)} )-A(\varrho (s),\mathcal{L}_{\varrho (s)}), u^{k_{L}}(s)-\varrho (s) \rangle _{B}   \nonumber
\\&~~~+2\langle f(u^{k_{L}}(s),\mathcal{L}_{u^{k_{L}}(s)} )-f(\varrho (s),\mathcal{L}_{\varrho (s)}), u^{k_{L}}(s)-\varrho (s) \rangle _{U_{1}}
\\&~~~+\left \| g(u^{k_{L}}(s),\mathcal{L}_{u^{k_{L}}(s)} )-g(\varrho (s),\mathcal{L}_{\varrho (s)}) \right \|_{\mathscr{L}(U_{2},U_{1})}^{2}
- \eta\left \| u^{k_{L}}(s)-\varrho (s) \right \|_{U_{1}}^{2} \nonumber
\\&~~~+2_{B ^{\ast}} \langle A(\varrho (s),\mathcal{L}_{\varrho (s)} ), u^{k_{L}}(s) \rangle _{B}+2\langle f(\varrho (s),\mathcal{L}_{\varrho (s)} ), u^{k_{L}}(s) \rangle _{U_{1}}   \nonumber
\\&~~~+2_{B ^{\ast}} \langle A(u^{k_{L}}(s),\mathcal{L}_{u^{k_{L}}(s)} )-A(\varrho (s),\mathcal{L}_{\varrho (s)} ), \varrho (s) \rangle _{B}
+2\langle f(u^{k_{L}}(s),\mathcal{L}_{u^{k_{L}}(s)} )-f(\varrho (s),\mathcal{L}_{\varrho (s)} ), \varrho (s) \rangle _{U_{1}}  \nonumber
\\&~~~+2\langle g(u^{k_{L}}(s),\mathcal{L}_{u^{k_{L}}(s)} ),g(\varrho (s),\mathcal{L}_{\varrho(s)} ) \rangle _{\mathscr{L}(U_{2},U_{1})}-\left \|g(\varrho (s),\mathcal{L}_{\varrho (s)} )  \right \|_{\mathscr{L}(U_{2},U_{1})}^{2} \nonumber
\\&~~~-2\eta\langle u^{k_{L}}(s), \varrho (s) \rangle _{U_{1}}
+\eta \left \|\varrho (s)  \right \|_{U_{1}}^{2} ]\text{d}s.   \nonumber
	\end{align}
  We can and will assume without loss of generality that there exists a sufficiently large constant $M$ such that $\mathbb{E} \left \| \varrho (z) \right \|^{2}_{U_{1}}\le M$.

By \textbf{(H4)}, \eqref{r9} and Jensen's Inequality, we have
\begin{align*}
		  \Pi (t):&=\mathbb{E} \int_{0}^{t}e^{-\eta  s}[2_{B ^{\ast}} \langle A(u^{k_{L}}(s),\mathcal{L}_{u^{k_{L}}(s)} )-A(\varrho (s),\mathcal{L}_{\varrho (s)}), u^{k_{L}}(s)-\varrho (s) \rangle _{B}   \nonumber
\\&~~~+2\langle f(u^{k_{L}}(s),\mathcal{L}_{u^{k_{L}}(s)} )-f(\varrho (s),\mathcal{L}_{\varrho (s)}), u^{k_{L}}(s)-\varrho (s) \rangle _{U_{1}}
\\&~~~+\left \| g(u^{k_{L}}(s),\mathcal{L}_{u^{k_{L}}(s)} )-g(\varrho (s),\mathcal{L}_{\varrho (s)}) \right \|_{\mathscr{L}(U_{2},U_{1})}^{2}
- \eta\left \| u^{k_{L}}(s)-\varrho (s) \right \|_{U_{1}}^{2}]\text{d}s
\\&\le\mathbb{E} \int_{0}^{t}e^{-\eta  s}[\lambda_{1}[\mathcal{G}(\left \| u^{k_{L}}(s)-\varrho (s) \right \| _{U_{1}}^{2})+d_{U_{1}}^{2} (\mathcal{L}_{u^{k_{L}}(s)} ,\mathcal{L}_{\varrho (s)})
- \eta\left \| u^{k_{L}}(s)-\varrho (s) \right \|_{U_{1}}^{2}]\text{d}s
\\&\le\left | \lambda _{1} \right |\mathcal{G}(\mathbb{E} \int_{0}^{t}\left \| u^{k_{L}}(s)-\varrho (s) \right \| _{U_{1}}^{2}\text{d}s).
	\end{align*}
Given any nonnegative function $\rho  \in L^{\infty }([0,T],\mathbb{R} )$ and letting $L \to \infty $, it follows from \eqref{p9} that
\begin{align}\label{p10}
		 &\mathbb{E}\int_{0}^{T}\rho (t) [e^{-\eta  t}\left \| u(t) \right \|_{U_{1}}^{2}  -\left \|x  \right \|_{U_{1}}^{2} ]\text{d}t  \nonumber
\\&\le\mathbb{E} \int_{0}^{T}\rho (t)\int_{0}^{t}e^{-\eta  s}[2_{B ^{\ast}} \langle A(\varrho (s),\mathcal{L}_{\varrho (s)} ), u(s) \rangle _{B}+2\langle f(\varrho (s),\mathcal{L}_{\varrho (s)} ), u(s) \rangle _{U_{1}} \nonumber
\\&~~~+2_{B ^{\ast}} \langle A^{*}(s)-A(\varrho (s),\mathcal{L}_{\varrho (s)} ), \varrho (s) \rangle _{B}
+2\langle f^{*}(s )-f(\varrho (s),\mathcal{L}_{\varrho (s)} ), \varrho (s) \rangle _{U_{1}}
\\&~~~+2\langle g^{*}(s ), g(\varrho (s),\mathcal{L}_{\varrho (s)} ) \rangle _{\mathscr{L}(U_{2},U_{1})}-\left \|g(\varrho (s),\mathcal{L}_{\varrho (s)} )  \right \|_{\mathscr{L}(U_{2},U_{1})}^{2}  \nonumber
\\&~~~-2\eta\langle u(s), \varrho (s) \rangle _{U_{1}}  +\eta \left \|\varrho (s)  \right \|_{U_{1}}^{2} ]\text{d}s\text{d}t+\left | \lambda _{1} \right |\int_{0}^{T}\rho (t)\mathcal{G}(\mathbb{E} \int_{0}^{t}\left \| u(s)-\varrho (s) \right \| _{U_{1}}^{2}\text{d}s)\text{d}t.   \nonumber
	\end{align}
Applying It\^o's formula to $e^{-\eta  t}\left \| u(t) \right \|_{U_{1}}^{2}  -\left \|x  \right \|_{U_{1}}^{2}$ implies
\begin{align}\label{p11}
		 &\mathbb{E}e^{-\eta  t}\left \| u(t) \right \|_{U_{1}}^{2}  -\left \|x \right \|_{U_{1}}^{2}   \nonumber
\\&=\mathbb{E} \int_{0}^{t}e^{-\eta  s}[2_{B ^{\ast}} \langle A^{*}(s ), u(s) \rangle _{B}+2\langle f^{*}(s ), u(s) \rangle _{U_{1}}
\\&~~~+\left \| g^{*}(s ) \right \|_{\mathscr{L}(U_{2},U_{1})}^{2}- \eta\left \| u(s) \right \|_{U_{1}}^{2}]\text{d}s.  \nonumber
	\end{align}
Substituting \eqref{p11} into \eqref{p10} gives
\begin{align}\label{p12}
		 &0\ge\mathbb{E} \int_{0}^{T}\rho (t)\int_{0}^{t}e^{-\eta  s}[2_{B ^{\ast}} \langle A^{*}(s )-A(\varrho (s),\mathcal{L}_{\varrho (s)} ), u(s)-\varrho (s) \rangle _{B} \nonumber
\\&~~~+2\langle  f^{*}(s )-f(\varrho (s),\mathcal{L}_{\varrho (s)} ), u(s)-\varrho (s) \rangle _{U_{1}}+\left \|g^{*}(s )-g(\varrho (s),\mathcal{L}_{\varrho (s)} )  \right \|_{\mathscr{L}(U_{2},U_{1})}^{2}
\\&~~~-\eta \left \|u(s)-\varrho (s)  \right \|_{U_{1}}^{2} ]\text{d}s\text{d}t-\left | \lambda _{1} \right |\int_{0}^{T}\rho (t)\mathcal{G}(\mathbb{E} \int_{0}^{t}\left \| u(s)-\varrho (s) \right \| _{U_{1}}^{2}\text{d}s)\text{d}t.   \nonumber
	\end{align}
Let $\varrho (s)=u(s)$, which implies $g^{*}(s )=g(\varrho _{s},\mathcal{L}_{\varrho _{s}} )$, $\text{d}t\times \mathbb{P}-\text{a.e.}$ by \eqref{p12}. Then, taking $\varrho =u-\xi\rho ^{*}v $ where $\rho ^{*} \in L^{\infty }([0 ,T],\mathbb{R} )$, $\xi>0$ and $v\in B$, we have
\begin{align*}
		 &0\ge\mathbb{E} \int_{0}^{T}\rho (t)\int_{0}^{t}e^{-\eta  s}[2_{B ^{\ast}} \langle A^{*}(s )-A(u(s)-\xi\rho ^{*}(s)v,\mathcal{L}_{u(s)-\xi\rho ^{*}(s)v} ), \xi\rho ^{*}(s)v \rangle _{B} \nonumber
\\&~~~+2\langle  f^{*}(s )-f(u(s)-\xi\rho ^{*}(s)v,\mathcal{L}_{u(s)-\xi\rho ^{*}(s)v} ), \xi\rho ^{*}(s)v \rangle _{U_{1}}-\eta \left \|\xi\rho ^{*}(s)v  \right \|_{U_{1}}^{2} ]\text{d}s\text{d}t
\\&~~~-\left | \lambda _{1} \right |\int_{0}^{T}\rho (t)\mathcal{G}(\xi^{2}\mathbb{E} \int_{0}^{t}\left \| \rho ^{*}v  \right \| _{U_{1}}^{2}\text{d}s)\text{d}t,  \nonumber
	\end{align*}
which implies
\begin{align*}
		 &\mathbb{E} \int_{0}^{T}\rho (t)\int_{0}^{t}e^{-\eta  s}[2_{B ^{\ast}} \langle A^{*}(s )-A(u(s)-\xi\rho ^{*}(s)v,\mathcal{L}_{u(s)-\xi\rho ^{*}(s)v} ), \rho ^{*}(s)v \rangle _{B} \nonumber
\\&~~~+2\langle  f^{*}(s )-f(u(s)-\xi\rho ^{*}(s)v,\mathcal{L}_{u(s)-\xi\rho ^{*}(s)v} ), \rho ^{*}(s)v \rangle _{U_{1}}-\eta \xi\left \|\rho ^{*}(s)v  \right \|_{U_{1}}^{2} ]\text{d}s\text{d}t
\\&~~~-\left | \lambda _{1} \right |\int_{0}^{T}\rho (t)\frac{\mathcal{G}(\xi^{2}\mathbb{E} \int_{0}^{t}\left \| \rho ^{*}v  \right \| _{U_{1}}^{2}\text{d}s)}{\xi}\text{d}t\le0.
	\end{align*}
  According to \textbf{(H1)}, \textbf{(H3)}, \eqref{r9} and Lebesgue's dominated convergence theorem,  and  letting $\xi\to 0$, we obtain
\begin{align*}
		 &\mathbb{E} \int_{0}^{T}\rho (t)\int_{0}^{t}e^{-\eta  s}[2_{B ^{\ast}} \langle A^{*}(s )-A(u(s),\mathcal{L}_{u(s)} ), \rho ^{*}(s)v \rangle _{B} \nonumber
\\&~~~+2\langle  f^{*}(s )-f(u(s),\mathcal{L}_{u(s)} ), \rho ^{*}(s)v \rangle _{U_{1}} ]\text{d}s\text{d}t\le0.
	\end{align*}
Similarly, the converse follows by letting $ \rho ^{*}(s)=- \rho ^{*}(s)$, and finally we can get
\begin{align*}
		 &\mathbb{E} \int_{0}^{T}\rho (t)\int_{0}^{t}e^{-\eta  s}[2_{B ^{\ast}} \langle A^{*}(s )-A(u(s),\mathcal{L}_{u(s)} ), \rho ^{*}(s)v \rangle _{B} \nonumber
\\&~~~+2\langle  f^{*}(s )-f(u(s),\mathcal{L}_{u(s)} ), \rho ^{*}(s)v \rangle _{U_{1}} ]\text{d}s\text{d}t=0.
	\end{align*}
 Then, the arbitrariness of $\rho ^{*}(s)$ and $v$ leads to
\begin{align}\label{p13}
A^{*}(s)=A(u(s  ),\mathcal{L}_{u(s  )}), \quad f^{*}(s)=f(u(s ) ,\mathcal{L}_{u(s )}), \quad \text{d}t\times \mathbb{P}-\text{a.e.}
	\end{align}
This completes the existence proof, i.e.,
\begin{align*}
u(t)=x+\int_{0}^{t}A(u(s  ),\mathcal{L}_{u(s  )})\text{d}s+\int_{0}^{t}f(u(s),\mathcal{L}_{u(s  )})\text{d}s+\int_{0}^{t}g(u(s  ),\mathcal{L}_{u(s  )})\text{d}W(s), \quad \text{d}t\times \mathbb{P}-\text{a.e.}
	\end{align*}

 The uniqueness of \eqref{r1} follows from the It$\hat{\text{o}} $ formula, $\textbf{(H4)}$ and \textbf{step 2} of Theorem 3.2. This completes the proof. \quad $\Box $
 ~\\
\\\textbf{Remark 3.3.}
It is noteworthy that condition \textbf{(H4)} falls within the category of local weak monotonicity conditions. While this condition imposes a less stringent regularity requirement on spatial variables compared to classical Lipschitz continuity, it nevertheless exhibits substantially greater regularity than the H$\ddot{\text{o}}$lder continuity condition. Given this observation, further relaxing the regularity requirements of SDE coefficients and systematically investigating the well-posedness of system \eqref{r1} under H$\ddot{\text{o}}$lder continuous coefficients is of profound theoretical importance. Within the classical Yamada-Watanabe framework, if the diffusion and drift coefficients of the system merely satisfy the H$\ddot{\text{o}}$lder continuity condition, the pathwise uniqueness of the weak solution generally cannot be established via the conventional contraction mapping principle or the martingale problem approach.  This fundamental challenge arises from the local oscillatory behavior of the solution under H$\ddot{\text{o}}$lder continuity, which may induce the phenomenon of non-measurable branching in the path space. Motivated by the generalized coupling method proposed in \cite{ref50}, further investigations into McKean-Vlasov SPDEs have been conducted to establish the weak uniqueness of the weak solution under H$\ddot{\text{o}}$lder continuity.
~\\
\\\textbf{(H4')} (H$\ddot{\text{o}}$lder continuity) Let $\alpha \in(0,1] $, $\beta \in(\frac{1}{2},1 ]$ and the embedding $B\subset U_{1}$ be compact. The map $A$ satisfies, for all $x,y\in B$, $\mu,\nu \in \mathcal{P}^{*}(U_{1})$,
\begin{align*}
		2_{B ^{\ast}} \langle A(x,\mu)-A(y,\nu), x-y\rangle _{B} \le \lambda_{1}(\left \| x-y \right \| _{U_{1}}^{2}+d_{U_{1}}^{2} (\mu ,\nu)),
	\end{align*}
and the functions $f$, $g$ satisfy,  for all $x,y\in U_{1}  $ and  $\mu, \nu\in \mathcal{P}^{*}(U_{1}  )$ with $\left \| x-y\right \| _{U_{1}}\vee d_{U_{1}} (\mu ,\nu)\le 1$,
\begin{align*}
		\left \langle f(x,\mu)-f(y,\nu),x-y \right \rangle _{U_{1}} \le\lambda_{1}(\left \| x-y \right \| _{U_{1}}^{\alpha+1 }+d_{U_{1}}^{2} (\mu ,\nu)) ,
	\end{align*}
\begin{align*}
	 \left \| g(x,\mu )-g(y,\nu ) \right \| _{\mathscr{L}(U_{2},U_{1})} \le\lambda_{2}(\left \| x-y\right \| _{U_{1}}^{\beta }+d_{U_{1}} (\mu ,\nu)).
	\end{align*}
\textbf{(H5)} For each \( x \in U_{1} \), \( g \) has a right inverse \( g^{-1} \) on \(U_{1}\) in the sense that
\[
g [g^{-1} x] = x \quad \text{for any } x \in U_{1},
\]
with \( g^{-1}[U_{1} ]\subset U_{2} \) and satisfies,
\[
\sup_{x \in U_{1}} \left| g^{-1}(t, x) \right|_{L(U_{1}, U_{2})} < \infty.
\]
Beyond the aforementioned conditions, the following lemma is also required:
~\\
\\\textbf{Lemma 3.3.} (Lemma B.1, \cite{ref50}) Let $V(t) \geq 0$ be an It\^o process with
\begin{align*}
    dV(t) = \eta(t) dt + dM(t),
\end{align*}
where $M$ is a continuous local martingale with quadratic variation
\begin{align*}
    \langle M \rangle (t) = \int_0^t m(s) ds, \quad t \geq 0.
\end{align*}
Assume that for some constants $A \geq 0$, $B > 0$, $\lambda > 0$ and a random variable $\varsigma \geq 0$
\begin{align*}
    \eta(t) \leq -\lambda V(t) + A, \quad m(t) \leq B, \quad t \leq \varsigma.
\end{align*}
Assume also that $\varsigma \leq T$ for some constant $T > 0$.
Then for every $\delta \in (0, 1/2)$ there exist constants $C_1, C_2 > 0$, which depend only on $\delta$ and $T$, such that
\begin{align*}
    P \left( \sup_{t \leq \varsigma} \left( V(t) - e^{-\lambda t} V(0) \right) \geq A \lambda^{-1} + B^{1/2} \lambda^{-\delta} R \right) \leq C_1 e^{-C_2 R^2}, \quad R \geq 0.
\end{align*}
\\\textbf{Theorem 3.5.} \emph{Consider \eqref{r1}. Suppose that the assumptions \textbf{(H1)}$-$\textbf{(H3)}, \textbf{(H4')} and \textbf{(H5)} hold.  Then for any initial value $x \in U_{1} $, system \eqref{r1} has  a weak solution  $u(t)$, and the weak solution $u(t)$  is unique in law; that is, any two such solutions with the same initial value have the same law.
~\\
\\\textbf{proof:}} We employ the Galerkin method to construct approximating solutions and subsequently establish the tightness of their laws in a suitable function space, thereby proving the existence of weak solutions. Following the approach in Theorem 3.1, we derive the finite-dimensional approximate system \eqref{r2}. To facilitate the proof, we first assume the existence of a weak solution $u^{k}(t)$  to the finite-dimensional system \eqref{r2} under conditions \textbf{(H1)}$-$\textbf{(H3)}, \textbf{(H4')} and \textbf{(H5)}. Moreover, this weak solution $u^{k}(t)$  is unique in law.

The existence of weak solutions for system \eqref{r1} closely mirrors the result presented in Theorem 3.2. Consequently, by  \eqref{r2}-\eqref{p13}, it follows that any weak limit point  $u(t)$ for the sequence $u^{k}(t)$, as $k\to \infty$ is a weak solution to \eqref{r1}. It is important to emphasize that the function $\mathcal{G}$ needs to be replaced by function $\mathcal{Y}(r)=r^{\frac{\alpha+1}{2}}+r^{\beta}$ in the proof process. Therefore, the central aspect of Theorem 3.4 lies in proving the uniqueness in law of the weak solution  $u(t)$. This proof is focused on establishing the weak uniqueness of the solution to \eqref{r1}. To establish the weak uniqueness of the solution to system \eqref{r1}, we will employ the concept of generalized coupling. From Proposition 3.3 of \cite{ref999} as well as Condition \textbf{(H4')},  we can construct sequences of operators $f^{m}$ and $g^{m}$ which
satisfy (one-sided) Lipschitz condition in $x$. Specifically, for any $x,y\in U_{1}  $ and  $\mu, \nu\in \mathcal{P}^{*}(U_{1}  )$ with $\left \| x-y\right \| _{U_{1}}\vee d_{U_{1}} (\mu ,\nu)\le 1$,
\begin{align*}
		\left \langle f^{m}(x,\mu)-f^{m}(y,\nu),x-y \right \rangle _{U_{1}} \le\lambda_{1}(\left \| x-y \right \| _{U_{1}}^{2 }+d_{U_{1}}^{2} (\mu ,\nu)) ,
	\end{align*}
\begin{align*}
	 \left \| g^{m}(x,\mu )-g^{m}(y,\nu ) \right \| _{\mathscr{L}(U_{2},U_{1})} \le\lambda_{2}(\left \| x-y\right \| _{U_{1}}+d_{U_{1}} (\mu ,\nu)).
	\end{align*}
And  $(f^{m},g^{m})$ converge uniformly to $(f,g)$ on any compact subset of $U_{1}$. Given any compact subset  $U^{'}_{1} \subset U_{1}$,  denote
\[
\tau_{1}=\inf \{t\ge 0;u(t)\notin U^{'}_{1}\}, \quad ,
 \Theta_{1}=\sup_{x\in U^{'}_{1}}\sup_{\mu\in\mathcal{P}^{*}(U_{1}) }\| f^{m}(x,\mu)-f(x,\mu)\|_{U_{1}},
\]
\[
\Theta_{2}=\sup_{x\in U^{'}_{1}}\sup_{\mu\in\mathcal{P}^{*}(U_{1}) }\| g^{m}(x,\mu)-g(x,\mu)\|_{\mathscr{L}(U_{2},U_{1})},\quad \Theta_{*}=\max\{\sqrt[\alpha]{\Theta_{1}} ,\sqrt[\beta]{\Theta_{2}}  \}.
\]
Now consider the following three McKean-Vlasov SPDEs:
\begin{align}
    \begin{cases}
        \text{d}u(t) = \left(A(u(t), \mathcal{L}_{u(t)}) + f(u(t), \mathcal{L}_{u(t)})\right) \text{d}t + g(u(t), \mathcal{L}_{u(t)}) \text{d}W(t), \\
        u(0) =  x \in U_1,
    \end{cases} \label{r1r}
\end{align}

\begin{align}
    \begin{cases}
        \text{d}u^{m}(t) = \left(A(u^{m}(t), \mathcal{L}_{u^{m}(t)}) + f^{m}(u^{m}(t), \mathcal{L}_{u^{m}(t)})\right) \text{d}t + g^{m}(u^{m}(t), \mathcal{L}_{u^{m}(t)}) \text{d}W(t), \\
        u(0) =x \in U_1,
    \end{cases} \label{r2r}
\end{align}
and
\begin{align}
    \begin{cases}
        \text{d}u^{*}_{m}(t) = & \, \left(A(u^{*}_{m}(t), \mathcal{L}_{u^{*}_{m}(t)}) + f^{m}(u^{*}_{m}(t), \mathcal{L}_{u^{*}_{m}(t)})+ \zeta(u(t) - u^{*}_{m}(t))\mathcal{X} _{t\le\tau}\right) \text{d}t \\
          & \,+ g^{m}(u^{*}_{m}(t), \mathcal{L}_{u^{*}_{m}(t)}) \text{d}W(t)
       \text{d}t, \\
        u(0) = & \, x \in U_1,
    \end{cases} \label{r3r}
\end{align}
where $\tau=\inf\{t\ge 0; \|u(t)-u^{*}_{m}(t)\|_{U_{1}}>2\Theta_{*}\}$ is some stopping time adapted to the filtration $(\mathcal{F}_t)_{t\geq 0}$ and $\varsigma=\Theta_{*}^{\vartheta-1} > 0$ is a constant with  $\vartheta\in (0, \alpha\wedge(2\beta-1))$.
%Since $G_n(t, x)$ converges to $G(t, x)$ for each $x \in H$, and $K$ is compact and $G(t, x)$ is uniformly continuous in $x$, we deduce that $\lim_{n \to \infty} \Delta_n^K = 0$ since all the estimates are time independent.
%
%Now define two stopping times
%\[
%\tau_n^K := \inf \left\{ t \geq 0 : |u(t) - u_n(t)| \geq 2 \Delta_n^K \right\},
%\]
%and
%\[
%\theta_K := \inf \left\{ t \geq 0 : u(t) \notin K \right\}.
%\]
%
%In the rest of the proof, we consider
%\[
%\tau = \tau_n^K \wedge \theta_K, \quad \lambda = (\Delta_n^K)^{\gamma - 1},
%\]
%where $\gamma$ will be specified later. We may slightly enlarge $\Delta_n^K$, but still assume
%$\lim_{n \to \infty} \Delta_n^K = 0$, such that $\Delta_n^K$ is much larger than the length scale for $G_n$ to be uniformly $\beta$-Hölder continuous in $x$, as discussed at the end of Section 3.2.

Now we estimate $u(t) - u_{m}^{*}(t)$ via pathwise arguments.  By Itô's formula, we obtain
\begin{align}\label{r5r}
&\left\| u(t) - u_{m}^{*}(t) \right\|_{U_{1}}^2\nonumber
\\&=  \int_0^t [2_{B ^{\ast}} \langle A(u(s),\mathcal{L}_{u(s)})-A(u_{m}^{*}(s),\mathcal{L}_{u_{m}^{*}(s)}), u(s) - u_{m}^{*}(s)\rangle _{B}\nonumber
\\&~~~+ 2\langle f(u(s),\mathcal{L}_{u(s)}) - f^{m}(u_{m}^{*}(s),\mathcal{L}_{u_{m}^{*}(s)}), u(s) - u_{m}^{*}(s) \rangle_{U_{1}}\nonumber
\\&~~~+ \left\| g(u(s),\mathcal{L}_{u(s)}) - g^{m}(u_{m}^{*}(s),\mathcal{L}_{u_{m}^{*}(s)}) \right\|_{\mathscr{L}(U_{2},U_{1})}^2- 2 \varsigma(u(s) - u^{*}_{m}(s))\mathcal{X} _{t\le\tau}]\text{d}s\nonumber
\\&~~~ + 2\int _{0}^{t}\langle u(s) - u^{*}_{m}(s),[g(u(s),\mathcal{L}_{u(s)}) - g^{m}(u_{m}^{*}(s),\mathcal{L}_{u_{m}^{*}(s)})]\text{d}W(s)\rangle_{U_{1}}\nonumber
\\&\le  \int_0^t [2_{B ^{\ast}} \langle A(u(s),\mathcal{L}_{u(s)})-A(u_{m}^{*}(s),\mathcal{L}_{u_{m}^{*}(s)}), u(s) - u_{m}^{*}(s)\rangle _{B}
\\&~~~+ 2\langle f(u(s),\mathcal{L}_{u(s)}) - f^{m}(u(s),\mathcal{L}_{u(s)}), u(s) - u_{m}^{*}(s) \rangle_{U_{1}}\nonumber
\\&~~~+ \left\| g(u(s),\mathcal{L}_{u(s)}) - g^{m}(u(s),\mathcal{L}_{u(s)}) \right\|_{\mathscr{L}(U_{2},U_{1})}^2\nonumber
\\&~~~+ 2\langle f^{m}(u(s),\mathcal{L}_{u(s)}) - f^{m}(u_{m}^{*}(s),\mathcal{L}_{u_{m}^{*}(s)}), u(s) - u_{m}^{*}(s) \rangle_{U_{1}}\nonumber
\\&~~~+\left\| g^{m}(u(s),\mathcal{L}_{u(s)})  - g^{m}(u_{m}^{*}(s),\mathcal{L}_{u_{m}^{*}(s)}) \right\|_{\mathscr{L}(U_{2},U_{1})}^2- 2 \varsigma\|u(s) - u^{*}_{m}(s)\|^{2}_{U_{1}}\mathcal{X} _{t\le\tau}]\text{d}s\nonumber
\\&~~~ + 2\int _{0}^{t}\langle u(s) - u^{*}_{m}(s),[g(u(s),\mathcal{L}_{u(s)}) - g^{m}(u(s),\mathcal{L}_{u(s)})]\text{d}W(s)\rangle_{U_{1}}\nonumber
\\&~~~ + 2\int _{0}^{t}\langle u(s) - u^{*}_{m}(s),[g^{m}(u(s),\mathcal{L}_{u(s)}) - g^{m}(u_{m}^{*}(s),\mathcal{L}_{u_{m}^{*}(s)})]\text{d}W(s)\rangle_{U_{1}}.\nonumber
\end{align}
Denote \( \tau^{*}= \tau \wedge \tau_1 \), then for \( t \leq \tau^{*} \) we have $u(t),u^{*}_{m}(t)\in U^{'}_{1}$ and
\[
\left\| u(t) - u^{*}_{m}(t) \right\|_{U_{1}} \leq 2\Theta_{*},
\]
which imply
\begin{align}\label{r6r}
&2_{B ^{\ast}} \langle A(u(s),\mathcal{L}_{u(s)})-A(u(s),\mathcal{L}_{u(s)}), u(s) - u_{m}^{*}(s)\rangle _{B}\nonumber
\\&+ 2\langle f(u(s),\mathcal{L}_{u(s)}) - f^{m}(u(s),\mathcal{L}_{u(s)}), u(s) - u_{m}^{*}(s) \rangle_{U_{1}}\nonumber
\\&+ \left\| g(u(s),\mathcal{L}_{u(s)}) - g^{m}(u(s),\mathcal{L}_{u(s)}) \right\|_{\mathscr{L}(U_{2},U_{1})}^2
\\&~~~\le C[\Theta_{*}^{2}+\Theta_{*}^{\alpha+1}+\Theta_{*}^{2\beta}].\nonumber
\end{align}
Then, by the conditions \textbf{(H4')} and \eqref{r9}, we get
\begin{align}\label{r7r}
&2\langle f^{m}(u(s),\mathcal{L}_{u(s)}) - f^{m}(u_{m}^{*}(s),\mathcal{L}_{u_{m}^{*}(s)}), u(s) - u_{m}^{*}(s) \rangle_{U_{1}}\nonumber
\\&~~~+\left\| g^{m}(u(s),\mathcal{L}_{u(s)})  - g^{m}(u_{m}^{*}(s),\mathcal{L}_{u_{m}^{*}(s)}) \right\|_{\mathscr{L}(U_{2},U_{1})}^2- 2 \varsigma\|u(s) - u^{*}_{m}(s)\|^{2}_{U_{1}}\mathcal{X} _{t\le\tau}]\nonumber
\\&\le C[\Theta_{*}^{\alpha}\|u(s) - u_{m}^{*}(s)\|_{U_{1}}+\Theta_{*}^{\alpha}\|u(s) - u_{m}^{*}(s)\|_{U_{1}}+\Theta_{*}^{2\beta}]-2 \varsigma\|u(s) - u^{*}_{m}(s)\|^{2}_{U_{1}}
\\&~~~+Cd ^{2}_{U_{1}} (\mathcal{L}_{u(s)},\mathcal{L}_{u_{m}^{*}(s)} )\nonumber
\\&\le C[\Theta_{*}^{\alpha}\|u(s) - u_{m}^{*}(s)\|_{U_{1}}+\Theta_{*}^{2\beta}+\Theta_{*}^{2}]-2 \varsigma\|u(s) - u^{*}_{m}(s)\|^{2}_{U_{1}}.\nonumber
\end{align}
By Young's inequality,
\begin{align*}
\Theta_{*}^{\alpha}\|u(s) - u_{m}^{*}(s)\|_{U_{1}}
\le \frac{1}{2}[\Theta_{*}^{\alpha+1}+\Theta_{*}^{\alpha-1}\|u(s) - u_{m}^{*}(s)\|_{U_{1}}].
\end{align*}
On the other hand, since \( \vartheta < \alpha \), there exists \( \Theta > 0 \) such that for any $\Theta_{*}\in (0,\Theta]$,
\begin{align}\label{r8r}
\Theta_{*}^{\alpha-1}-2\Theta_{*}^{\vartheta-1}\le-\Theta_{*}^{\vartheta-1}.
\end{align}
Further, substituting \eqref{r6r} and \eqref{r8r} into \eqref{r5r}, we obtain that for any \( t \in[0, \tau^{*}] \) and $\Theta_{*}\in (0,\Theta]$,
\begin{align}\label{r9r}
&\text{d}[\left\| u(t) - u_{m}^{*}(t) \right\|_{U_{1}}^2]\nonumber
\\&\le C[-\Theta_{*}^{\vartheta-1}\|u(t) - u_{m}^{*}(t)\|^{2}_{U_{1}}
+\Theta_{*}^{\alpha+1}+\Theta_{*}^{2\beta}+\Theta_{*}^{2}]\text{d}t+ C\Theta_{*}^{\beta+1}\text{d}W(s).
\end{align}
Now we apply Lemma 3.4 with a fixed \( T > 0 \), then we obtain
\[
 V(t) = \left\| u(t) - u_{m}^{*}(t) \right\|_{U_{1}}^2, \quad \zeta = \tau^{*} \wedge T,\quad \lambda = \Theta_{*}^{\vartheta-1},
\]
\[
 A = C \left(\Theta_{*}^{\alpha+1}+\Theta_{*}^{2\beta}+\Theta_{*}^{2} \right), \quad B = C \Theta_{*}^{2\beta+2},
\]
and
\[
\mathbb{P}\left( \sup_{t \leq \zeta}( \left\|u(t) - u_{m}^{*}(t) \right\|_{U_{1}}^2 )  \geq A\lambda^{-1}+B^{\frac{1}{2}}\lambda^{-\delta}R \right)\leq C_{1}e^{-C_{2}R^{2}},\quad R\geq 0,
\]
for  every $\delta\in (0,\frac{1}{2})$.

Recall that \( \vartheta < \alpha \) and \( \vartheta < 2\beta - 1 \), hence there exists \( \ell> 0 \) such that for \( \Theta > 0 \) small enough, we obtain
\begin{align*}
[\Theta_{*}^{\alpha+1}+\Theta_{*}^{2\beta}+\Theta_{*}^{2}]\Theta_{*}^{1-\vartheta} =[\Theta_{*}^{\alpha+2-\vartheta}+\Theta_{*}^{2\beta+1-\vartheta}+\Theta_{*}^{3-\vartheta}]
\le\varepsilon\Theta_{*}^{2+2\ell}.
\end{align*}
Further we know that there exists $\delta_{0}\in (0,\frac{1}{2})$ such that $\beta+\delta_{0}(1-\vartheta)>1$. And let $\ell\in (0,\frac{\beta+\delta_{0}(1-\vartheta)-1}{3})$, then we have
\[
B^{\frac{1}{2}}\lambda^{-\delta} \leq \varepsilon \Theta_{*}^{2 + 3\ell}, \quad \Theta_{*} \in [0, \Theta].
\]
Let  \( R = \Theta_{*}^{-\ell} \),  and denote

\[
\mathcal{A} = \left\{ \omega : \sup_{t \leq \tau^{*} \wedge T}   \left\|u(t) - u_{m}^{*}(t) \right\|_{U_{1}}  \geq\Theta_{*}^{1 + \ell} \right\}.
\]
Then,  we conclude that
\[
\mathbb{P}(\mathcal{A}) \leq C_1 e^{-C_2 \Theta_{*}^{-2\ell} }.
\]
On the other hand, taking \( \Theta \leq 1 \), we have for \( \Theta_{*} \in (0,  \Theta] \), $\Theta_{*}^{\ell}\leq 1$. Then, on the set \( \Omega\backslash\mathcal{A} \), we have
\[
 \left\|u(t) - u_{m}^{*}(t) \right\|^{2}_{U_{1}}  \leq\Theta_{*}^{2 + 2\ell}\leq4\Theta_{*}^{2}, \quad t \leq \tau^{*} \wedge T.
\]
Recall that the stopping time is defined as \( \tau=\inf\{t\ge 0; \|u(t)-u^{*}_{m}(t)\|_{U_{1}}>2\Theta_{*}\} \). Given that  \( u(t),u^{*}_{m}(t) \) exhibit continuous trajectories in the function space, it follows that on \( \Omega\backslash\mathcal{A} \), the stopping time satisfies \( \tau > \tau^{*} \wedge T \). Consequently, we establish that
\(
\tau^{*} \wedge T = \tau_{1} \wedge T.
\)
We thus derive the following result:
\[
\mathbb{P }\left(\sup_{t \in [0, \tau_{1}\wedge T]}   \left\|u(t) - u_{m}^{*}(t) \right\|_{U_{1}}  \geq\Theta_{*}^{1 + \ell} \right) \leq C_1 e^{-C_2 \Theta_{*}^{-2\ell} }.
\]
This implies, for any $\varepsilon > 0$, as $m$ gets large and $\Theta_{*} \to 0$, we have
\begin{align}\label{r11r}
P \left( \sup_{t \in [0, \tau_{1}\wedge T]} \left\|u(t) - u_{m}^{*}(t) \right\|_{U_{1}}  > \varepsilon \right) \to 0, \quad m \to \infty.
\end{align}

In the following, we can employ Girsanov's theorem to estimate $u^{m}(t)- u_{m}^{*}(t)$. Let  $d_{TV}$ represent the total variation distance on $C([0,T]; U_{1})$, and define $KL(\cdot \| \cdot)$ as the Kullback-Leibler divergence between two probability measures on $C([0,T]; U_{1})$.  Since the coefficients  \( f^{m} \) and \( g^{m} \) satisfy strong monotonicity condition, it follows that  \( u^{m}(t) \) and \( u_{m}^{*}(t) \) are the strong solutions to \eqref{r2r} and \eqref{r3r}, respectively. Consequently, they can be interpreted as the images of \( W^{*} \), \( W \) under a measurable mapping. Specifically,   the SPDE solved by $u_{m}^{*}(t)$ can be derived from  the SPDE solved by $u^{m}(t)$ by introducing an additional drift term into the Brownian sample path:
\[
\text{d}W^{*}(t) = \text{d}W(t) + \mathbb{G} (t) \text{d}t, \quad \mathbb{G} (t) = g^{m}(u^{m}(t))^{-1} \varsigma(u(t) - u^{*}_{m}(t))\mathcal{X} _{t\le\tau\wedge T} \text{d}t.
\]
This gives
\[
d_{\text{TV}} \left( \text{Law}(u^{m}(t)| _{[0,T]}), \text{Law}(u_{m}^{*}(t) | _{[0,T]}) \right) \leq d_{\text{TV}} \left( \text{Law}(W(t) | _{[0,T]}), \text{Law}(W^{*} (t)| _{[0,T]}) \right),
\]

Since $(g)^{-1}$ is a bounded operator, recalling the definition of the stopping time $\tau$ and   control term $\varsigma$, we deduce that
\[
\mathbb{E} \left[ \int_0^\infty \left\| \mathbb{G}(t) \right\|_{U_{1}}^2 \, dt \right] \leq C_{T} \Theta_{*}^{2\varsigma}.
\]
By Pinsker's inequality,
\[
KL (\text{Law}(W^{m}(t)| _{[0,T]}) \| \text{Law}(W^{*}(t)|_{ [0,T]}))\ge \frac{1}{2}[d_{TV}(\text{Law}(W^{m}(t) | _{[0,T]}), \text{Law}(W^{*}(t)| _{[0,T]}))]^{2} .
\]Then, the following bound for the total variation distance holds:
\begin{align}\label{r10r}
d_{\text{TV}} \left( \text{Law}(u^{m}(t)| _{[0,T]}), \text{Law}(u_{m}^{*}(t) | _{[0,T]}) \right)  \leq C _{T} \Theta_{*}^{\varsigma}.
\end{align}

We can now finish the proof of uniqueness in law. Consider $\mathcal{W}$ a real-valued, bounded continuous function defined on $C([0, T], U_{1})$. By \eqref{r10r}, we have
\[
\mathbb{E}\mathcal{W}(u^{m}(t)| _{[0,T]}) - \mathbb{E}\mathcal{W}(u_{m}^{*}(t)|_{[0,T]}) \to 0, \quad m \to \infty.
\]
In addition, \eqref{r11r} implies that on $\{ \tau_{1} \geq T \} \cap (\Omega \setminus \mathcal{A})$, we have
 \[u_{m}^{*}(t)| _{[0,T]} \to u(t) | _{[0,T]},\quad  m \to \infty,
 \]
 in probability in $C([0, T], U_{1})$.
Since $\lim_{n \to \infty} P(\mathcal{A}) = 0$, we conclude that
\[
\limsup_{m \to \infty} |\mathbb{E}\mathcal{W}(u_{m}^{*}(t)| _{[0,T]}) - \mathbb{E}\mathcal{W}(u(t)|_{[0,T]})| \leq 2 \sup_x |\mathcal{W}(x)| \mathbb{P}(\tau_{1}\leq T),
\]
which implies
\begin{align}\label{r13r}
\limsup_{m \to \infty} |\mathbb{E}\mathcal{W}(u^{m}(t)| _{[0,T]}) - \mathbb{E}\mathcal{W}(u(t)|_{[0,T]})| \leq 2 \sup_x |\mathcal{W}(x)| \mathbb{P}(\tau_{1} \leq T).
\end{align}
Similar to \eqref{p1}-\eqref{p5}, we obtain that for any $p > 1$, there exists a constant $M>0$ such that
\begin{equation}\label{eq:condition}
\mathbb{E}\left[ \sup_{t \in [0,T]} \|u(t)\|_{U_{1}}^2 \right] + \mathbb{E}\left[ \int_0^T \|u(s)\|_{B}^{p} \, ds \right] \leq M .
\end{equation}
For any \( \varepsilon > 0 \), let \( R = \left(\frac{2M}{\varepsilon}\right)^{1/p} \). By Chebyshev's inequality:
\[
\mathbb{P}\left( \int_0^T \|u(s)\|_{B}^p ds > R^p \right) \leq \frac{\mathbb{E}\left[ \int_0^T \|u(s)\|_{B}^p ds \right]}{R^p} \leq \frac{M}{R^p} = \frac{\varepsilon}{2},
\]
which implies
\begin{equation}\label{eq:integral}
\mathbb{P}\left( \int_0^T \|u(s)\|_{B}^p ds \leq R^p \right) \geq 1 - \frac{\varepsilon}{2}.
\end{equation}
Given that \( \mathbb{E}[\sup_{t \in [0,T]} \|u(t)\|_{U_{1}}^2] \leq M \), assuming condition \textbf{(H3)} holds,  it follows as a standard procedure, based on the Kolmogorov continuity theorem, that there exists a parameter  \( \varpi \in (0,1/2) \) and a random variable \( \Re \) satisfying \( \mathbb{E}\Re^2 \leq  M \) such that for any $t,s\in[0,T]$,
\[
\|u(t) \|_{U_{1}}\le \Re,\quad \|u(t) - u(s)\|_{U_{1}} \leq \Re|t-s|^\varpi \quad \text{a.s.}
\]
Applying Chebyshev to \( \Re \):
\[
\mathbb{P}\left( \Re > K \right) \leq \frac{\mathbb{E}\Re^2}{K^2} \leq \frac{ M}{K^2}, \quad \text{for } K > 0.
\]
Choosing \( K = \sqrt{\frac{2M}{\varepsilon}} \), we have
\begin{equation}\label{eq:holder}
\mathbb{P}\left( \Re \leq K \right) \geq 1 - \frac{\varepsilon}{2}.
\end{equation}
Then define the space:
\[
\mathcal{S} = L^p([0,T], B) \cap C^{0,\varpi}([0,T], U_{1}).
\]
By the Aubin-Lions lemma and compact embedding \( B \subset U_{1} \), bounded sets in \( \mathcal{S}  \) are compact in \( C([0,T], U_{1}) \). Define
\[
\mathcal{E} = \left\{ u \in C([0,T], H) \,\Big|\, \|u\|_{C^{0,\varpi}} \leq K, \, \int_0^T \|u(s)\|_{B}^p ds \leq R^p \right\},
\]
then \( \mathcal{E} \) is compact in \( C([0,T], U_{1}) \). For any \( t \in [0,T] \), consider the projection map \( \pi_t: C([0,T], U_{1}) \to U_{1} \), \( \pi_t(u) = u(t) \). Since \( \mathcal{E} \) is compact and \( \pi_t \) is continuous, the image
\[
\mathcal{U} = \pi_t(\mathcal{E}) = \{ u(t) \in U_{1} \mid u \in \mathcal{E} \}
\]
is compact in \( U_{1} \). Combining \eqref{eq:integral} and \eqref{eq:holder}:
\[
\mathbb{P}\left( u \in \mathcal{E} \right) \geq 1 - \frac{\varepsilon}{2} - \frac{\varepsilon}{2} = 1 - \varepsilon.
\]
Therefore
\begin{align}\label{r15r}
\mathbb{P}\left( u(t) \in \mathcal{U} \right) \geq 1 - \varepsilon, \quad \forall t \in [0,T],
\end{align}
which implies that for any $\varepsilon>0$, there exists a compact subset $\mathcal{U}\subset U_{1}$ such that with probability at least $1-\varepsilon$, $u(t)$ never leaves $\mathcal{U}$ for all $t\in[0, T]$. Consequently,  we can find a compact set $U^{'}_{1}$ such that $\mathbb{P}(\tau_{1} \leq T)$ is arbitrarily
small. Letting $\varepsilon\to 0$ and applying \eqref{r13r} and \eqref{r15r},
we finally deduce that
\[\mathbb{E}\mathcal{W}(u^{m}(t)| _{[0,T]}) - \mathbb{E}\mathcal{W}(u(t)|_{[0,T]})\to 0,\quad m\to \infty,\]
which completes the proof of weak uniqueness. This follows from the fact that an arbitrary weak solution  to \eqref{r1r} is
now uniquely specified on any finite time interval $[0,T]$ as the weak limit of the solutions to \eqref{r2r}, and hence the law of $u$ is uniquely determined.
~\\
\\\textbf{Remark 3.6.} The initial step in establishing Theorem 3.5 involves addressing the existence and uniqueness of weak solutions for the finite-dimensional system \eqref{r2}, derived via Galerkin projection techniques. In fact, under conditions \textbf{(H1)}$-$\textbf{(H3)}, the existence of a weak solution to system \eqref{r2} can be established via the first step of Theorem 3.1. The proof of the uniqueness of the weak solution follows the same reasoning as that in Theorem 3.5. Specifically, by employing a similar argument as in \eqref{r1r}-\eqref{r15r}, we establish that the weak solution to \eqref{r2} is unique in the sense of law. Since the proof is relatively straightforward, we omit the details for brevity.

\section{\textup{Ergodicity and exponential mixing under Lyapunov conditions}}
In this section, we explore  the long-term behavior of the solution of system \eqref{r1}, which encompass aspects such as existence, uniqueness, exponential convergence, and exponential mixing of invariant measures. The assumption of dissipation is frequently crucial when examining the pertinent properties of the invariant measures. Specifically, it is necessary to assume that the coefficients $f$ and $g$  are Lipschitz continuous. To achieve ergodicity and exponential mixing of \eqref{r1} using It$\hat{\text{o}}$'s formula, a stronger assumption needs to be made:
\begin{align*}
		\left \langle f(x,\mu)-f(y,\nu),x-y \right \rangle _{U_{1}} \le-\lambda_{2}(\left \| x-y \right \| _{U_{1}}^{2}+d_{U_{1}}^{2} (\mu ,\nu)) ,
	\end{align*}
where $\lambda_{2}$ is sufficiently large compared to the Lipschitz constant for the diffusion coefficient. However, this  imposes a very stringent structural restriction. Therefore,  we further improve the assumptions in the following way. We demonstrate that when the coefficients satisfy weak
monotonicity assumptions, we can utilize a Lyapunov function to analyze invariant measure properties. We find this approach to be straightforward and intriguing in its own right. Specifically, we attach dissipative conditions to the Lyapunov function rather than to the structure of the system itself. This expansion significantly broadens the applicability of our approach. Throughout this paper  Section 4 and  Section 5, we will assume that \textbf{(H1)}$-$\textbf{(H5)} hold or \textbf{(H1)}$-$\textbf{(H3)}, \textbf{(H4')} and \textbf{(H5)} hold.

 Denote by $C^{2}(U_{1} ,\mathbb{R}^{+} )$ the family of all real-valued  nondecreasing functions $V(x )$ defined on $U_{1}$, which are twice continuously differentiable and $V(0)=0$. Applying It$\hat{\text{o}} $ formula, it follows for $u(t)$, a solution to \eqref{r1}, that
\begin{align*}
V(u(t))&=V(x)+\int_{0}^{t}[ _{B ^{\ast}} \langle A(u(s), \mathcal{L}_{u(s)}), \nabla V(u(s))\rangle _{B}+\left \langle f(u_{s}, \mathcal{L}_{u_{s}}), \nabla V(u(s)) \right \rangle _{U_{1}}
\\&~~~+\frac{1}{2} \nabla ^{2}V(u(s))\left \| g(u_{s}, \mathcal{L}_{u_{s}}) \right \|^{2}_{\mathscr{L}(U_{2},U_{1})}]\text{d}s+\int_{0}^{t}  \left \langle \nabla V(u(s)),g(u_{s}, \mathcal{L}_{u_{s}})\text{d}W(s) \right \rangle _{U_{1}}.
	\end{align*}
We define an operator $LV:U_{1}\to \mathbb{R} $  associated with \eqref{r1} as follows: for any $x \in U_{1} $
\begin{align*}
LV(x)=_{B ^{\ast}} \langle A(x, \mathcal{L}_{x}), \nabla V(x)\rangle _{B}+\left \langle f(x, \mathcal{L}_{x}), \nabla V(x) \right \rangle _{U_{1}}
+\frac{1}{2} \nabla ^{2}V(x)\left \| g(x, \mathcal{L}_{x}) \right \|^{2}_{\mathscr{L}(U_{2},U_{1})},
	\end{align*}
then
\begin{align*}
V(u(t))&=V(x)+\int_{0}^{t}LV(u(s))\text{d}s+M_{t},
	\end{align*}
where $M_{t}=\int_{0}^{t}  \left \langle \nabla V(u(s)),g(u(s), \mathcal{L}_{u(s)})\text{d}W(s) \right \rangle _{U_{1}}$ is a continuous martingale with $M_{0}=0$.

Let $\mathcal{P}_{V} $ be the family of all probability measures on $(U_{1}, \mathcal{B}(U_{1}))$ and satisfy $\mu (V)<\infty $ for any $\mu \in \mathcal{P}(U_{1})$. Then the following Wasserstein quasi-distance in the metric space is then induced by the Lyapunov function $V$:
\begin{align*}
\Pi _{V}(\mu ,\nu  )=\underset{\pi \in C(\mu ,\nu )}{\inf} \int _{U_{1} \times U_{1}}V(x -y )\pi (\text{d}x ,\text{d}y ),
	\end{align*}
where $C(\mu ,\nu )$ denotes the set of all couplings between $\mu$ and $\nu$, see \cite{ref3} for more details on this Wasserstein quasi-distance. According to the Kantorovich–Rubinshtein theorem, we know that $\Pi _{V}(\mu ,\nu  )$ has the following alternative expression:
\begin{align*}
		\Pi _{V}(\mu ,\nu  ):=\sup\{\left | \int  \Phi \text{d}\mu- \int  \Phi  \text{d}\nu \right |:\left \|  \Phi   \right \|_{LipV}\le 1  \},
	\end{align*}
 where $ \Phi :U_{1}\to \mathbb{R} $ and $\left \|  \Phi  \right \|_{LipV}:=\sup_{x\ne y}\frac{| \Phi (x)- \Phi (y) |}{V(x-y)} $.

We denote by  $$p(t, 0, x, \Gamma  ) :=\mathbb{P}  (\omega :u(t; 0, x)\in\Gamma )$$
 the transition probability for the  solution $u(t;0,x)$, where $u(0;0,x)=x$ and $\Gamma\in \mathcal{B}(U_{1})$. For any $t\ge0$ we associate a mapping $P^{*}_{t}:\mathcal{P}(U_{1} )\to  \mathcal{P}(U_{1} )$ defined by
\begin{align}\label{p21}
		P^{*}_{t}\mu (\Gamma )=\hat{p}(t,0)\mu(\Gamma )=\int _{U_{1} }p(t,0,x ,\Gamma )\mu(\text{d}x ).
	\end{align}
For any $F\in C_{b}(U_{1})$, which is defined as the the set of all bounded continuous functions $F : U_{1}\to \mathbb{R}$ endowed with the norm $\left \| F \right \| _{\infty }=\sup_{u \in U_{1} } \left \| F(u ) \right \| $, we define the following semi-group $P_{0,t}$ for $t\ge 0$
\begin{align}\label{p22}
		P_{0,t}F(x )=\int _{U_{1} }F(y )p(t,0,x ,\text{d}y  ).
	\end{align}
In particular, the  operator $P_{0,t}$ is written as $P_{t}$. Firstly, we will establish some properties for the solution $u(t;0,x)$ of \eqref{r1}.
~\\
 \\\textbf{Theorem 4.1.} \emph{Consider \eqref{r1}. The following statement holds: the solution $u(t;0,x)$ is a time homogeneous Markov process in $U_{1}$ and  has the Feller property.
 ~\\
\\\textbf{proof:}}  Let us divide this proof into two steps:
\\\textbf{step 1:} (Markov and Feller property) Note that the Markov and the  Feller property for $u(t;0,x)$ is obtained from the same property for $u^{k_{L}}(t;0,x)$ by the usual approximation argument. Similar to Theorem 2.1 in \cite{ref1}, we obtain that $u^{k_{L}}(t;0,x)$ is a Markov process and has the Feller property. Hence for any bounded Borel measurable function $h : U_{1}\to \mathbb{R} $ and $0 \le s\le  t <\infty $,
\begin{align*}
		\mathbb{E} (h(u^{k_{L}}(t))\mid \mathcal{F} _{s}^{k_{L}})=\mathbb{E} (h(u^{k_{L}}(t))\mid u^{k_{L}}(s)),
	\end{align*}
where $\mathcal{F} _{t}^{k_{L}}=\sigma  \{W^{k_{L}}(s):0\le s\le t \}$. Let $L\to \infty $, then we get that  $u^{k_{L}}(t)\to u(t)$ in $\mathfrak{L} ^{2}([0,T]\times \Omega ,U_{1})$ and  $W^{k_{L}}(t)\to W(t)$, which implies
\begin{align*}
		\mathbb{E} (h(u(t))\mid \mathcal{F} _{s})=\mathbb{E} (h(u(t))\mid u(s)).
	\end{align*}
This proves that $u(t)$ is a  Markov process.

For any $F\in C_{b}(U_{1})$, from the definition of the semi-group $P_{t}$, we have $\left \| P_{t}F \right \| _{\infty}\le\left \| F \right \| _{\infty }<\infty $. The next major task is to prove that $ P_{t}F$ is continuous. We just need to prove that for any sequence $x_{n}\in U_{1} $, $x \in  U_{1} $, when $\lim_{n \to \infty } \left \| x _{n}-x   \right \| _{ U_{1} }=0$, we have $\lim_{n \to \infty }  | P_{t}F(x  _{n})- P_{t}F(x  )   |=0$.

Let $P^{k_{L}}_{t}F(x )=\int _{U_{1} }F(y )p^{k_{L}}(t,0,x ,\text{d}y  )$ and $p^{k_{L}}(t, 0, x, \Gamma  ) :=\mathbb{P}  (\omega :u^{k_{L}}(t; 0, x^{k_{L}})\in\Gamma )$. Since $u^{k_{L}}(t)\to u(t)$ in $\mathfrak{L} ^{2}([0,T]\times \Omega ,U_{1})$, we have for any $x\in U_{1}$
\begin{align*}
		P^{k_{L}}_{t}F(x )\to P_{t}F(x ), \quad L\to \infty.
	\end{align*}
Hence
\begin{align*}
		&\lim_{n \to \infty } | P_{t}F(x  _{n})- P_{t}F(x  )   |\\&\le\lim_{n \to \infty } [| P^{k_{L}}_{t}F(x  _{n})- P^{k_{L}}_{t}F(x  )   |+ | P_{t}F(x  _{n})- P^{k_{L}}_{t}F(x _{n}  )   |+| P_{t}F(x )- P^{k_{L}}_{t}F(x  )   |]\\&\le\lim_{n \to \infty } [| P_{t}F(x  _{n})- P^{k_{L}}_{t}F(x _{n}  )   |+| P_{t}F(x )- P^{k_{L}}_{t}F(x  )   |].
	\end{align*}
Let  $L\to \infty$, then we finally obtain $P_{t}F\in C_{b}(U_{1})$.
~\\
\\\textbf{step 2:} (Time-homogeneous) According to the definition of Time-homogeneous, We just need to prove that for all $t\ge0$, $\theta\ge0$, $x\in U_{1}$ and $ \Gamma \in \mathcal{B} (U_{1})$,
\begin{align*}
		p(t ,0 ,x,\Gamma)=p(t+\theta ,\theta ,x,\Gamma).
	\end{align*}
Firstly,
\begin{align}\label{p31}
		u(t;0,x)=x+\int_{0}^{t}A(u(s  ),\mathcal{L}_{u(s  )})\text{d}s+\int_{0}^{t}f(u(s  ),\mathcal{L}_{u(s  )})\text{d}s+\int_{0}^{t}g(u(s  ),\mathcal{L}_{u(s  )})\text{d}W(s),
	\end{align}
in addition, $u(t+\theta ;\theta,x)$ is determined by the solution $u(t)$
\begin{align}\label{p32}
		u(t+\theta ;\theta,x)&=x+\int_{\theta }^{t+\theta }A(u(s  ),\mathcal{L}_{u(s  )})\text{d}s+\int_{\theta}^{t+\theta}f(u(s  ),\mathcal{L}_{u(s  )})\text{d}s
\\&~~~+\int_{\theta }^{t+\theta }g(u(s  ),\mathcal{L}_{u(s  )})\text{d}W(s),\nonumber
	\end{align}
This equation is equivalent to
\begin{align}\label{p333}
		u(t+\theta ;\theta,x)&=x+\int_{0}^{t }A(u(s+\theta   ),\mathcal{L}_{u(s+\theta   )})\text{d}s+\int_{0}^{t}f(u(s+\theta   ),\mathcal{L}_{u(s+\theta   )})\text{d}s \nonumber
\\&~~~+\int_{0 }^{t }g(u(s+\theta   ),\mathcal{L}_{u(s+\theta   )})\text{d}\tilde{W} (s),
	\end{align}
where $\tilde{W} (t)=W(t+\theta)-W(\theta)$ is a  cylindrical Wiener process with the same distribution as $W(t)$. Then we see by \eqref{p31}, \eqref{p32}, \eqref{p333} and  the uniqueness of \eqref{r1},
\begin{align*}
		\mathbb{P}  (\omega :u(t; 0, x)\in\Gamma )=\mathbb{P}  (\omega :u(t+\theta; \theta, x)\in\Gamma ),
	\end{align*}
i.e.,
\begin{align*}
		p(t ,0 ,x,\Gamma)=p(t+\theta ,\theta ,x,\Gamma).
	\end{align*}
This completes the proof. $\Box $
~\\

 In order to analyze the existence, uniqueness, exponential convergence, and exponential mixing of invariant measures, we need the following conditions $\textbf{(H6)}$ and $\textbf{(H7)}$, which are used in many works.
~\\
\\\textbf{(H6)} (Dissipative Lyapunov condition) There exists constant $\delta >0$ such that for any $x,y\in U_{1}$, $\mu,\nu \in \mathcal{P}^{*}(U_{1})$ and $\pi \in C(\mu ,\nu )$
\begin{align}\label{p33}
		\int _{U_{1}\times U_{1}}\widehat{L} V(x-y)\pi(\text{d}x,\text{d}y)\le-\delta \int _{U_{1}\times U_{1}}V(x-y)\pi(\text{d}x,\text{d}y),
	\end{align}
where
\begin{align*}
\widehat{L} V(x-y)&=_{B ^{\ast}} \langle A(x, \mu)-A(y, \nu), \nabla V(x-y)\rangle _{B}+\left \langle f(x, \mu)-f(y, \nu), \nabla V(x-y) \right \rangle _{U_{1}}
\\&~~~+\frac{1}{2} \nabla ^{2}V(x-y)\left \| g(x,\mu)-g(y,\nu) \right \|^{2}_{\mathscr{L}(U_{2},U_{1})}.
	\end{align*}
%\textbf{(H4*)} (The following stronger version of\textbf{(H4)}) There exist constants $\lambda^{*}>0 $ and $\beta \in(\frac{1}{2},1 ]$.  For all $x,y\in B$, $\mu,\nu \in \mathcal{P}^{*}(U_{1})$,
%\begin{align*}
%		2_{B ^{\ast}} \langle A(u,\mu)-A(v,\nu), u-v\rangle _{B} \le -\lambda^{*}(\left \| u-v \right \| _{U_{1}}^{2}+d_{U_{1}}^{2} (\mu ,\nu)) ,
%	\end{align*}
%and   the functions $f$, $g$ satisfy,  for all $x,y\in U_{1}  $ and  $\mu, \nu\in \mathcal{P}^{*}(U_{1} )$ with $\left \| x-y \right \| _{U_{1}}\vee\left \| \mu-\nu  \right \| _{U_{1}}\le 1$,
%\begin{align*}
%		\left \langle f(x,\mu)-f(y,\nu),x-y \right \rangle _{U_{1}} \le\lambda_{1}(\left \| x-y \right \| _{U_{1}}^{2\beta }+d_{U_{1}}^{2\beta} (\mu ,\nu)) ,
%	\end{align*}
%\begin{align*}
%	 \left \| g(x,\mu )-g(y,\nu ) \right \| _{\mathscr{L}(U_{2},U_{1})} \le\lambda_{1}(\left \| x-y \right \| _{U_{1}}^{\beta }+d_{U_{1}}^{\beta} (\mu ,\nu)).
%	\end{align*}
\\\textbf{(H7)}  There exist constants $\delta^{*}, \delta_{0}, c >0$ such that for any $x\in U_{1}$ and  $\mu\in \mathcal{P}^{*}(U_{1})$
\begin{align}\label{p33}
		\int _{U_{1}}L V(x)\mu(\text{d}x)\le-\delta^{*} \int _{U_{1}}V(x)\mu(\text{d}x)+\delta_{0},
	\end{align}
and
\begin{align*}
V(x-y)\le M[V(cx)+V(cy)].
	\end{align*}

We have the following results on the existence and uniqueness, exponential convergence, and exponential mixing of invariant measures for \eqref{r1}.
~\\
\\\textbf{Theorem 4.2.}
 \emph{
	\begin{enumerate}[(\textbf{\uppercase\expandafter{\romannumeral1}})]
		\item Assume that \textbf{(H6)} holds, then
\begin{enumerate}[(\textbf{a})]
		\item For any initial value $x,y\in U_{1}$ such that $\mu_{0}=\mathcal{L} _{x}\in \mathcal{P}_{V}$ and  $\nu_{0}=\mathcal{L} _{y}\in \mathcal{P}_{V}$, we have that for any $t\ge0$
\begin{align*}
\Pi _{V}(P^{*}_{t}\mu_{0},P^{*}_{t}\nu_{0})\le e^{-\delta t }\Pi _{V}(\mu_{0},\nu_{0});
	\end{align*}
\end{enumerate}
\begin{enumerate}[(\textbf{b})]
        \item there exists $\nu ^{*}\in\mathcal{P}_{V} $ such that
        \begin{align*}
\sup\{\Pi _{V}(P^{*}_{t}\nu ^{*},\nu ^{*}):t\ge0\}<\infty .
	\end{align*}
And there exists a unique measure $\mu^{*}\in\mathcal{P}_{V}$  satisfying
\begin{align*}
		P^{*}_{t}\mu^{*}(\Gamma )=\mu^{*}(\Gamma ),
	\end{align*}
for any $t\ge0$ and $\Gamma \in \mathcal{B}(U_{1} )$, i.e., $\mu^{*}$ is a invariant measure and satisfies
\begin{align*}
\Pi _{V}(P^{*}_{t}\nu,\mu^{*})\le e^{-\delta t }\Pi _{V}(\nu,\mu^{*}),
	\end{align*}
for any $\nu\in\mathcal{P}_{V}$.
	\end{enumerate}
	\end{enumerate}
\begin{enumerate}[(\textbf{\uppercase\expandafter{\romannumeral2}})]
		\item Under assumptions $\textbf{(H6)}$-$\textbf{(H7)}$, the  invariant measure $\mu^{*}$ of \eqref{r1} is uniformly exponential mixing in the sense of Wasserstein metric. More precisely, for any $t_{0}>0$, $t\ge t_{0}$ and $\nu \in \mathcal{P} _{V}$,
\begin{align*}
		\Pi _{V}(P^{*}_{t}\nu,\mu^{*})\le Me^{-\delta t } [\frac{\delta_{0}}{\delta ^{*}(1-e^{- t_{0}\delta ^{*}})}+\int_{U_{1} }V( cx   )\nu(\text{d}x )].
	\end{align*}
	\end{enumerate}
~\\
\textbf{Proof of (\textbf{\uppercase\expandafter{\romannumeral1}}):}} Similar to the proof of Theorem 4.1 in \cite{ref3}. This argument is again quite standard, we thus omit the details.
~\\
 \emph{\\\textbf{Proof of (\textbf{\uppercase\expandafter{\romannumeral2}}):}}
 For all $F\in C_{b}(\mathcal{R})$,
\begin{align}\label{r111}
		P_{t}F(x  )&=\int _{U_{1} }F(y )p(t,0,x   ,\text{d}y  )\nonumber
\\&=\int _{U_{1} }F(y )\mathbb{P} (\omega :u(t;0,x  )\in\text{d}y )
\\&=\mathbb{E} F(u(t;0,x  )),\nonumber
	\end{align}
and by the proof of Theorem 4.1 in \cite{ref3}, we have
\begin{align}\label{r112}
		\mathbb{E} V( u(t;0,x )-u(t;0,y ) )\le e^{-\delta t }\mathbb{E} V( x-y  ).
	\end{align}
 By \eqref{r111}, \eqref{r112}, Chapman–Kolmogorov equation and $\textbf{(H7)}$, we obtain that for any $\nu \in \mathcal{P}_{V}$,
\begin{align}\label{g10}
		&\Pi _{V}(P^{*}_{t}\nu,\mu^{*}) \nonumber
\\&=\Pi _{V}(P^{*}_{t}\nu,P^{*}_{t}\mu^{*})\nonumber
\\&=\sup_{\left \| F \right \|_{LipV}\le1 } \left \| \int_{U_{1} }  F(z)\int_{U_{1} }p(t,0,x, \text{d}z )\nu(\text{d}x )- \int_{U_{1} } F(z)\int_{U_{1} }p(t,0,y, \text{d}z )\mu ^{*}(\text{d}y ) \right \|\nonumber
\\&=\sup_{\left \| F \right \|_{LipV}\le1 } \left \|  \int_{U_{1} }\mathbb{E} F(u(t;0,x ))\nu(\text{d}x )- \int _{U_{1} }\mathbb{E} F(u(t;0,y )) \mu ^{*}(\text{d}y )\right \|
\\&\le\sup_{\left \| F \right \|_{LipV}\le1 }\{\left \| F \right \|_{LipV}(  \int_{U_{1} }\int_{U_{1} }\mathbb{E} V( u(t;0,x )-u(t;0,y ) )\mu ^{*}(\text{d}y)\nu(\text{d}x ) )\}\nonumber
\\&\le e^{-\delta t } \int_{U_{1} }\int_{U_{1} }V( x -y  )\mu ^{*}(\text{d}y)\nu(\text{d}x )\nonumber
\\&\le Me^{-\delta t } [\int_{U_{1} }V( cx   )\nu(\text{d}x )+\int_{U_{1} }V( cy  )\mu ^{*}(\text{d}y)]\nonumber .
	\end{align}
 Applying It$\hat{\text{o}} $ formula to $V(cu(t;0,x))$ and by $\textbf{(H7)}$, there exists a constant $c^{*}$
 \begin{align*}
c^{*}=\left\{\begin{matrix}
  c^{2}&  c\ge1, \\
  c&  c<1,
\end{matrix}\right.
	\end{align*}
such that
 \begin{align*}
\mathbb{E} V(cu(t;0,x))&=V(cx)+\mathbb{E} \int_{0}^{t}LV(cu(s;0,x))\text{d}s
\\&\le V(cx)-c ^{*}\delta ^{*} \mathbb{E} \int_{0}^{t}V(cu(s;0,x))\text{d}s+\delta_{0}t.
	\end{align*}
Let
\begin{align*}
		v (t)= V(cx)-c ^{*}\delta ^{*}  \int_{0}^{t}v (s)\text{d}s+\delta_{0}t,
	\end{align*}
which implies that $v(t)$ satisfies the following equation
\begin{align*}
		\dot{v }(t) =\delta_{0}-c ^{*}\delta ^{*}v (t),
	\end{align*}
with initial condition $v (0)=V(cx)$. Solving this equation for $v (t)$, we get
\begin{align*}
		v (t)=V(cx)e^{-c ^{*}\delta ^{*}t}+\frac{\delta_{0}}{c ^{*}\delta ^{*}}(1-e^{-c ^{*}\delta ^{*}t}) .
	\end{align*}
Hence we obtain $ \mathbb{E} V(cu(t;0,x))\le V(cx)e^{-c^{*}\delta^{*}t}+\frac{\delta_{0}}{c ^{*}\delta ^{*}}$ by comparison principle. Then
\begin{align*}
		\int_{U_{1} }V( cy  )\mu ^{*}(\text{d}y)&=\int_{U_{1} }V( cy  )P^{*}_{t}\mu ^{*}(\text{d}y)
\\&=\int_{U_{1} }V(cy)\int_{U_{1} }p(t,0,z, \text{d}y )\mu ^{*}(\text{d}z )
\\&=\int_{U_{1} }\mathbb{E} V(cu(t;0,z))\mu ^{*}(\text{d}z )
\\&\le e^{-c^{*}\delta^{*}t}\int_{U_{1} }V( cz  )\mu ^{*}(\text{d}z)+\frac{\delta_{0}}{c ^{*}\delta ^{*}}.
	\end{align*}
For any $t_{0}>0$, we obtain
\begin{align*}
		\int_{U_{1} }V( cy  )\mu ^{*}(\text{d}y)\le \frac{\delta_{0}}{c^{*}\delta^{*}(1-e^{- t_{0}c ^{*}\delta ^{*}})},
	\end{align*}
which implies
\begin{align*}
		\Pi _{V}(P^{*}_{t}\nu,\mu^{*}) \le Me^{-\delta t } [\frac{\delta_{0}}{c ^{*}\delta ^{*}(1-e^{- t_{0}c ^{*}\delta ^{*}})}+\int_{U_{1} }V( cx   )\nu(\text{d}x )].
	\end{align*}
$\Box$

~\\
%\\\textbf{Remark 4.8.} Assume the entrance measure $\mu _{t}$ with respect to $P_{r}^{*}(r\ge0)$ is periodic with periodic $T$ or quasi-periodic with periodic $T_{1}$ and $T_{2}$, where the reciprocals of $T_{1}$ and $T_{2}$ are rationally linearly independent. Then the invariant measures are $$\mu=\frac{1}{T} \int_{0}^{T} \mu _{t}\text{d}t \quad \text{and} \quad \mu =\frac{1}{T_{1}T_{2}} \int_{0}^{T_{1}}\int_{0}^{T_{2}} \hat{\mu }_{t_{1},t_{2}}\text{d}t_{1}\text{d}t_{2}.$$
\section{\textup{Limit theorems of McKean-Vlasov SPDEs}}
In this section, based on the  uniform exponential mixing of the  invariant measure $\mu^{*}$, we further obtain SLLN and CLT of McKean-Vlasov SPDEs. Before presenting  some details, to facilitate presentation to follow, we introduce the following preliminaries. Let us fix a weight function $\mathcal{H}(r)>0$ for any $r\ge0$, which is increasing, continuous and bounded, and let $C_{\mathcal{H} }(U_{1})$ denote the family of all continuous functionals on $U_{1}$ such that
\begin{align*}
\begin{split}
		\left \| \Phi \right \| _{\mathcal{H}  }:=\sup_{x\ne y}\frac{| \Phi (x)- \Phi (y) |}{V(x-y)\cdot [\mathcal{H}(\left \| x \right \|_{U_{1}})+\mathcal{H}(\left \| y \right \|_{U_{1}})]}+\sup_{x \in U_{1} }\frac{|\Phi(x)|}{\mathcal{H}(\left \| x  \right \|_{U_{1} } )  }<\infty.
\end{split}
	\end{align*}
 For any $\Phi\in C_{\mathcal{H} }(U_{1})$, we set
$$\Psi_{t}^{x}[\Phi]=\int_{0}^{t}\Phi(u(s;x))\text{d}s, \quad \psi_{t}^{x}[\Phi]=\frac{1}{t}\Psi_{t}[\Phi],$$
where $u(t;x)=u(t;0,x)$ is the solution of \eqref{r1}.
\subsection{SLLN}
Firstly, based on the uniformly exponential mixing of the  measure $\mu^{*}$ of \eqref{r1}, we prove the SLLN for  a class of McKean-Vlasov SPDEs under weak monotonicity conditions. Assume, in addition, that the following stronger version of $V(x)$ holds true:
~\\
\\\textbf{(H8)}  There exist constants $\delta_{1}, \delta_{2}\in \mathbb{R} $ and $\kappa >0$ such that for any $x\in U_{1}$ and  $\mu\in \mathcal{P}^{*}(U_{1})$
\begin{align*}
		\int _{U_{1}} |\nabla  V(x)|\left \| g(x, \mu) \right \|_{\mathscr{L}(U_{2},U_{1})}\mu(\text{d}x)\le\delta_{1} \int _{U_{1}}V^{\kappa }(x)\mu(\text{d}x)+\delta_{2}.
	\end{align*}
~\\
\\\textbf{Theorem 5.1.}\emph{ Let assumptions \textbf{(H6)}$-$\textbf{(H7)} and \textbf{(H8)} with $\kappa =1$ hold. Assume also that for some $k\in \mathbb{N} ^{+}$ such that $\hat{\delta } >0$, where
$$\hat{\delta } =2k\delta^{*}-(2k-1)[\delta_{0}+2\delta ^{2}_{1}k+2\delta^{2}_{2}(k-1)].$$
  Then for any $x\in U_{1}$ and $\Phi\in C_{\mathcal{H} }(U_{1} ) $, we obtain the following conclusions:
\begin{enumerate}[(1)]
        \item There exists a constant $M>0$ such that
        \begin{align}\label{g16}
\begin{split}
		&\mathbb{E} \left |\Psi_{t}^{x}[\Phi] -(\Phi,\mu ^{*}) \right | ^{2k}
\\&\le 2k(2k-1)M^{2}\delta ^{-1}[\frac{\delta_{0}}{c ^{*}\delta ^{*}(1-e^{- c ^{*}\delta ^{*}})}+V( cx   )]^{2k}\left \| \Phi \right \|^{2k}_{\mathcal{H}}t^{-k}, \quad t\ge 1,
\end{split}
	\end{align}
where $(\Phi,\mu ^{*})=\int _{U_{1} }\Phi(z) \mu ^{*}(\text{d}z )$;
		\item  There exists a constant $M>0$ such that
\begin{align}\label{g17}
\begin{split}
		\left |\psi_{t}^{x}[\Phi] -(\Phi,\mu ^{*}) \right |\le M\left \| \Phi\right \|_{\mathcal{H}}t^{-\frac{1}{2} +\varepsilon }, \quad t\ge T_{\varepsilon}(\omega )\ge 1,\quad \mathbb{P} -a.s.,
\end{split}
	\end{align}
where  the random time $T_{\varepsilon}(\omega )$ is $\mathbb{P}-$a.s. finite. Moreover,
\begin{align}\label{g18}
		\mathbb{E} T_{\varepsilon}^{k }(\omega )&\le k(k+1)(2k(2k-1)M^{2}\delta ^{-1})^{k}M_{k}[(\frac{\delta_{0}}{c ^{*}\delta ^{*}(1-e^{- c ^{*}\delta ^{*}})})^{2k}\\&~~~+V^{2k}(cx )+\frac{[\delta_{0}+2\delta^{2}_{2}(2k-1)]}{c^{*}\hat{\delta } }]\left \| \Phi \right \|^{k(k+1)}_{\mathcal{H}} \nonumber.
	\end{align}
	\end{enumerate}
\textbf{proof}} For any given $\Phi\in C_{\mathcal{H} }(U_{1} ) $ and $x \in U_{1} $, it follows from Theorem 4.2 and Chapman-Kolmogorov align that
\begin{align}\label{g15}
\left | P_{t}\Phi(x )- (\Phi,\mu ^{*})\right | &=\left |\int _{U_{1} } \Phi(z )p(t,0,x ,\text{d}z )- \int _{U_{1} }\Phi(z )\mu ^{*}(\text{d}z )\right |\nonumber
\\&=\left | \int _{U_{1} }\int _{U_{1} }\Phi(z  )p(t,0,y  ,\text{d}z )p(0,0,x ,\text{d}y  )- \int _{U_{1} }\Phi(z )\mu ^{*}(\text{d}z )\right |\nonumber
\\&=\left | \int _{U_{1} }\Phi(z  )P^{*}_{t}p(0,0,x ,\text{d}z )- \int _{U_{1} }\Phi(z )\mu ^{*}(\text{d}z)\right |
\\&\le \left \| \Phi \right \| _{LipV}\cdot \Pi _{V}(P^{*}_{t}\mu _{0}^{x },\mu ^{*})\nonumber
\\&\le M\left \| \Phi \right \| _{\mathcal{H}}e^{-\delta t } [\frac{\delta_{0}}{c ^{*}\delta ^{*}(1-e^{- c ^{*}\delta ^{*}})}+\int_{U_{1} }V( cy   )p(0,0,x ,\text{d}y )]\nonumber
\\&= M\left \| \Phi \right \| _{\mathcal{H}}e^{-\delta t } [\frac{\delta_{0}}{c ^{*}\delta ^{*}(1-e^{- c ^{*}\delta ^{*}})}+V( cx   )]\nonumber,
	\end{align}
where $\mu _{0}^{x }(\Gamma )=p(0,0,x,\Gamma )$ for any $\Gamma \in \mathcal{B}(\mathcal{R} )$.  In addition, given  $k \ge 1$, in view of It\^o's formula, \textbf{(H7)}, \textbf{(H8)} and Young's inequality, we get
\begin{align*}
&\mathbb{E} (V^{2k}(u(t)))\\&=V^{2k}(x )+2k\mathbb{E} \int_{0}^{t}V^{2k-1}(u(s))[ _{B ^{\ast}} \langle A(u(s), \mathcal{L}_{u(s)}), \nabla V(u(s))\rangle _{B}+\left \langle f(u_{s}, \mathcal{L}_{u_{s}}), \nabla V(u(s)) \right \rangle _{U_{1}}
\\&~~~+\frac{1}{2} \nabla ^{2}V(u(s))\left \| g(u_{s}, \mathcal{L}_{u_{s}}) \right \|^{2}_{\mathscr{L}(U_{2},U_{1})}]\text{d}s\\&~~~+k(2k-1)\mathbb{E}\int_{0}^{t}V^{2k-2}(u(s))  |\nabla V(u(s))|^{2}\left \| g(u_{s}, \mathcal{L}_{u_{s}}) \right \|^{2}_{\mathscr{L}(U_{2},U_{1})}]\text{d}s
\\&\le V^{2k}(x )+[2\delta^{2}_{1}k(2k-1)-2k\delta ^{*}]\int_{0}^{t}V^{2k}(u(s))\text{d}s+2k\delta_{0}\int_{0}^{t}V^{2k-1}(u(s))\text{d}s
\\&~~~+2\delta^{2}_{2}k(2k-1)\int_{0}^{t}V^{2k-2}(u(s))\text{d}s
\\&\le V^{2k}(x )+[2\delta^{2} _{1}k(2k-1)-2k\delta^{*}+\delta_{0}(2k-1)+2\delta^{2}_{2}(2k-1)(k-1)]\int_{0}^{t}\mathbb{E}V^{2k}(u(s))\text{d}s
\\&~~~+[\delta_{0}+2\delta^{2}_{2}(2k-1)]t.
	\end{align*}
Then, it follows from Gronwall's lemma that
\begin{align}\label{i1}
\mathbb{E} (V^{2k}(u(t)))\le e^{-\hat{\delta } t}V^{2k}(x )+\frac{[\delta_{0}+2\delta^{2}_{2}(2k-1)]}{\hat{\delta } } <\infty.
	\end{align}
Hence  the assumptions in Lemma 2.1 of \cite{ref4} hold for $\mathbb{B}=U_{1}$, $\varphi (t)=Me^{-\delta t }$ and $\psi ( x  )=[\frac{\delta_{0}}{c ^{*}\delta ^{*}(1-e^{- c ^{*}\delta ^{*}})}+V( cx   )]$, then the desired assertion \eqref{g16} holds.

In addition, the assumptions in Definition 2.5 and Proposition 2.6 of \cite{ref5} hold for $\gamma (t)=Me^{-\delta t }$, $\rho (\left \| x \right \|_{U_{1}} )=\frac{\delta_{0}}{c ^{*}\delta ^{*}(1-e^{- c ^{*}\delta ^{*}})}+V( cx   )$,
and $$\sigma(\left \| x \right \|_{U_{1}})= M_{k}[(\frac{\delta_{0}}{c ^{*}\delta ^{*}(1-e^{- c ^{*}\delta ^{*}})})^{2k}+V^{2k}(cx )+\frac{[\delta_{0}+2\delta^{2}_{2}(2k-1)]}{c^{*}\hat{\delta } }],$$
then the desired conclusion (2) holds.
The proof of Theorem 5.1 is complete. $\quad \Box$
\subsection{Central limit theorem}
%To state the CLT, we first introduce the corrector $\Phi  :C_{\mathcal{K} }(\mathcal{R} )\longrightarrow \mathbb{R} $ such that
%\begin{align*}
%\Phi  [F(x )] =\int_{0}^{\infty } [P_{t}F(x )-\int _{\mathcal{R}}F(\phi )\mu ^{*}(\text{d}\phi ) ]\text{dt},
%	\end{align*}
In this subsection, let us fix $x\in U_{1} $ and an arbitrary function $\Phi\in C_{\mathcal{H} }(U_{1})$ such that $\int _{U_{1}}\Phi(x )\mu ^{*}(\text{d}x )=0$ and set
\begin{align*}
\mathcal{M} (t)^{x}_{\Phi}&=\int_{0}^{\infty } [\mathbb{E} (\Phi(u(s;x ))|\mathscr{F}_{t} )-\mathbb{E}(\Phi(u(s;x ))|\mathscr{F}_{0} ) ]\text{d}s
\\&=\int_{0}^{\infty } [\mathbb{E} (\Phi(u(s;x ))|\mathscr{F}_{t} )-P_{t}\Phi(x ) ]\text{d}s.
\end{align*}
For the property of $\mathcal{M} (t)^{x}_{\Phi}$, we have the following lemma:
~\\
\\\textbf{Lemma 5.2.} \emph{Assume the conditions of Theorem 5.1 hold. Then  $\mathcal{M} (t)^{x}_{\Phi}$ is a well-defined zero-mean martingale with $\mathbb{E} \left | \mathcal{M} (t)^{x}_{\Phi} \right |^{k} <\infty $.
\\\textbf{proof}} Firstly, by \eqref{g15} we have
\begin{align}\label{g20}
\begin{split}
		\left \| P_{t}\Phi(x ) \right \|  &\le M\left \| \Phi \right \| _{\mathcal{H}}e^{-\delta t } [\frac{\delta_{0}}{c ^{*}\delta ^{*}(1-e^{- c ^{*}\delta ^{*}})}+V( cx   )],
\end{split}
	\end{align}
then, by  the dominated convergence theorem,  for any $t>r\ge0$, we obtain
\begin{align*}
\begin{split}
\mathbb{E} (\mathcal{M} (t)^{x}_{\Phi}|\mathscr{F}_{r}  )&=\int_{0}^{\infty } [\mathbb{E}( \mathbb{E} (\Phi(u(s;x ))|\mathscr{F}_{t} )|\mathscr{F}_{r}) -\mathbb{E} (\mathbb{E}(\Phi(u(s;x ))|\mathscr{F}_{0} )|\mathscr{F}_{r})  ]\text{d}s\\
&=\int_{0}^{\infty } [\mathbb{E} (\Phi(u(s;x ))|\mathscr{F}_{r} )-\mathbb{E}(\Phi(u(s;x ))|\mathscr{F}_{0} ) ]\text{d}s\\
&=\mathcal{M} (r)^{x}_{\Phi}.
\end{split}
\end{align*}
In addition, it's not hard to get that $\mathbb{E} \mathcal{M} (t)^{x}_{\Phi}=0$.

Next, by \eqref{r111} and the  Markov properties of  $u(t)$, we obtain
\begin{align}\label{g21}
\begin{split}
\mathcal{M} (t)^{x}_{\Phi}&=\int_{0}^{\infty } [\mathbb{E} (\Phi(u(s;x ))|\mathscr{F}_{t} )-\mathbb{E}(\Phi(u(s;x ))|\mathscr{F}_{0} ) ]\text{d}s\\
&=\int_{0}^{t }\Phi(u(s;x )) \text{d}s+\int_{t}^{\infty } \mathbb{E} (\Phi(u(s;x ))|\mathscr{F}_{t} )\text{d}s-\int_{0}^{\infty }\mathbb{E}(\Phi(u(s;x )) ) \text{d}s\\
&=\int_{0}^{t }\Phi(u(s,x )) \text{d}s+\int_{t}^{\infty } P_{s-t}\Phi(u(s;x ))\text{d}s-\int_{0}^{\infty }\mathbb{E}(\Phi(u(s;x )) ) \text{d}s.
\end{split}
\end{align}
Thus, \eqref{i1}, \eqref{g20} and the definition of $\left \| \cdot  \right \| _{\mathcal{H} }$ yield
\begin{align*}
&\mathbb{E} \left | \mathcal{M} (t)^{x}_{\Phi} \right | ^{k}\\&
=\mathbb{E} \left |\int_{0}^{t }\Phi(u(s;x )) \text{d}s+\int_{t}^{\infty } P_{s-t}\Phi(u(s;x ))\text{d}s-\int_{0}^{\infty }P_{s}\Phi(x  ) \text{d}s\right | ^{k}\\
&=M_{k}\mathbb{E} \left | \left \| \Phi \right \|_{\mathcal{H} }[\int_{0}^{t}h(u(s;x ))\text{d}s
+M\int_{t}^{\infty }[\frac{\delta_{0}}{c ^{*}\delta ^{*}(1-e^{- c ^{*}\delta ^{*}})}+V( cu(s;x)   )]e^{-\delta (s-t) } \text{d}r]\right |^{k}
\\&~~~+M_{k}\mathbb{E} \left |\int_{0}^{\infty }M\left \| \Phi \right \|_{\mathcal{H} }[\frac{\delta_{0}}{c ^{*}\delta ^{*}(1-e^{- c ^{*}\delta ^{*}})}+V( cx   )]e^{-\delta s } \text{d}r) ]  \right |^{k}\\
&\le M(t) \left \| \Phi \right \|^{k}_{\mathcal{H} }[1+V^{k}(cx ) ],
\end{align*}
where $M(t)$ is an increasing function, which implies that $\mathbb{E} \left \| \mathcal{M} (t)^{x}_{\Phi} \right \| ^{k}< \infty  $ for $t\in \mathbb{R} $.\quad $\Box$
~\\

For any integer $j\ge 1$, by \eqref{g21}, we have
\begin{align}\label{g22}
\begin{aligned}
\mathcal{M} (j)^{x}_{\Phi}=\mathcal{M} (j-1)^{x}_{\Phi}+\int_{j-1}^{j} \Phi(u(s;x ))\text{d}s+\int_{0}^{\infty } P_{s}\Phi(u(j;x ))-P_{s}\Phi(u(j-1;x ))\text{d}s.
\end{aligned}
\end{align}
Next consider the conditional variance for $\{\mathcal{M} (n)^{x}_{\Phi} \}_{n\in\mathbb{Z}_{+}    }$,
\begin{align*}
\begin{aligned}
\Upsilon  ^{2}[ \mathcal{M} (n)^{x}_{\Phi}] =\sum_{j=1}^{n} \mathbb{E} ((\mathcal{M} (j)^{x}_{\Phi} -\mathcal{M} (j-1)^{x}_{\Phi})^{2}|\mathscr{F}_{j-1} ).
\end{aligned}
\end{align*}
By \eqref{g22} and the Markov properties of  $u(t)$, we obtain
\begin{align*}
\begin{aligned}
&\mathbb{E} ((\mathcal{M} (j)^{x}_{\Phi} -\mathcal{M} (j-1)^{x}_{\Phi})^{2}|\mathscr{F}_{i-1} )\\&=\mathbb{E} (\left \| \int_{i-1}^{i} \Phi(u(s;x ))\text{d}s+\int_{0}^{\infty } P_{s}\Phi(u(j;x ))-P_{s}\Phi(u(j-1;x ))\text{d}s\right \|^{2} |\mathscr{F}_{j-1}).
\end{aligned}
\end{align*}
Then let
\begin{align}\label{a5}
\Pi_{1} [\Phi(x )] =\int_{0}^{\infty } P_{t}\Phi(x )\text{dt},
\end{align}
and
 \begin{align}\label{a6}
 \Pi_{2} [\Phi(x )]=\mathbb{E} \left | \int_{0}^{1}\Phi(u(s;x ))\text{d}s+\Pi_{1}  [\Phi(u(1;x ))]-\Pi_{1}  [\Phi(x )]  \right |^{2},
  \end{align}
  which implies that by  Markov property
\begin{align}\label{g23}
\begin{aligned}
\Upsilon  ^{2}[ \mathcal{M} (n)^{x}_{\Phi}] =\sum_{j=1}^{n}\Pi_{2} [\Phi(u(j-1;x) )].
\end{aligned}
\end{align}

To establish the CLT, we  derive the following  crucial lemmas.
~\\
\\\textbf{Lemma 5.3.} \emph{Let assumptions \textbf{(H6)}$-$\textbf{(H7)} and \textbf{(H8)} with $\kappa =\frac{1}{2} $ hold. Then there exist constants $c_{1},c_{2}>0$ such that
\begin{align*}
\begin{aligned}
 \mathbb{E} (\sup_{z \in [0,t]}e^{c_{1}V(u(z;x) ) } )\le e^{c_{2}(1+V( x)  )},
\end{aligned}
\end{align*}
for any $t\ge1$ and $x \in U_{1} $.
\\\textbf{proof} }We first show $\mathbb{E} e^{c_{1}V(u(t;x))}<\infty $ for any $c_{1}>0$ and $t\ge 1$.  In view of It\^o's formula and $\textbf{(H7)}$,
\begin{align*}
 0\le V(u(t;x))&=V(x)+\int_{0}^{t}[ _{B ^{\ast}} \langle A(u(s), \mathcal{L}_{u(s)}), \nabla V(u(s))\rangle _{B}+\left \langle f(u_{s}, \mathcal{L}_{u_{s}}), \nabla V(u(s)) \right \rangle _{U_{1}}
\\&~~~+\frac{1}{2} \nabla ^{2}V(u(s))\left \| g(u_{s}, \mathcal{L}_{u_{s}}) \right \|^{2}_{\mathscr{L}(U_{2},U_{1})}]\text{d}s+\int_{0}^{t}  \left \langle \nabla V(u(s)),g(u_{s}, \mathcal{L}_{u_{s}})\text{d}W(s) \right \rangle _{U_{1}}
\\&\le V(x)-\delta ^{*}  \int_{0}^{t}V(u(s;x))\text{d}s+\delta_{0}t+ \widehat{\mathcal{M}}(t) ,
\end{align*}
where $ \widehat{\mathcal{M}}(t) =\int_{0}^{t}  \left \langle \nabla V(u(s)),g(u_{s}, \mathcal{L}_{u_{s}})\text{d}W(s) \right \rangle _{U_{1}}$. Then for some $\iota_{1} >0$ to be determined later, we have
\begin{align}\label{a1}
 \mathbb{E} e^{ \iota_{1} \int_{0}^{t}V(u(s;x))\text{d}s }&\le e^{\frac{\iota_{1} }{\delta ^{*}} (V(x)+\delta _{0}t)}\mathbb{E} e^{\frac{\iota_{1} }{\delta ^{*}}\int_{0}^{t}  \left \langle \nabla V(u(s)),g(u_{s}, \mathcal{L}_{u_{s}})\text{d}W(s) \right \rangle _{U_{1}}}\nonumber
 \\&\le e^{\frac{\iota_{1} }{\delta ^{*}} (V(x)+\delta _{0}t)}[\mathbb{E} e^{\frac{2\iota_{1}^{2} }{(\delta ^{*})^{2}}\int_{0}^{t} |\nabla  V(u(s))|^{2}\left \| g(u(s), \mathcal{L}_{u_{s}}) \right \|^{2}_{\mathscr{L}(U_{2},U_{1})}\text{d}s}]^{\frac{1}{2} }
 \\&\le e^{\frac{\iota_{1} }{\delta ^{*}} (V(x)+\delta _{0}t+\frac{2\iota_{1}\delta^{2} _{2} }{(\delta ^{*})}t)}[\mathbb{E} e^{\frac{4\iota_{1}^{2}\delta^{2} _{1} }{(\delta ^{*})^{2}}\int_{0}^{t} V(u(s;x))\text{d}s}]^{\frac{1}{2} }.\nonumber
\end{align}
Let $\iota_{1}=\frac{(\delta ^{*})^{2} }{4\delta^{2} _{1}}$, which implies that for any $t\ge 1$,
\begin{align}\label{a2}
 \mathbb{E} e^{ \iota_{1} \int_{0}^{t}V(u(s;x))\text{d}s }\le e^{\frac{2\iota_{1} }{\delta ^{*}} (V(x)+\delta _{0}t+\frac{2\iota_{1}\delta^{2} _{2} }{(\delta ^{*})}t)}<\infty.
\end{align}
In addition, let $\iota_{2}=\frac{\delta ^{*} }{4\delta _{1}}$, and we have by \eqref{a1} and \eqref{a2}
\begin{align*}
 \mathbb{E}  e^{\iota _{2}V(u(z;x))}&\le e^{\iota_{2} (V(x)+\delta _{0}t)}\mathbb{E}  e^{\iota_{2}\widehat{\mathcal{M}}(t)}
 \\&\le e^{\iota_{2} (V(x)+\delta _{0}t+2\iota_{2}\delta^{2} _{2}t)}[\mathbb{E} e^{4\iota_{2}^{2}\delta^{2} _{1}\int_{0}^{t} V(u(s;x))\text{d}s}]^{\frac{1}{2} }
 \\&\le e^{\iota_{2} (V(x)+\delta _{0}t+2\iota_{2}\delta^{2} _{2}t)+\frac{\iota_{1} }{\delta ^{*}} (V(x)+\delta _{0}t+\frac{2\iota_{1}\delta^{2} _{2} }{(\delta ^{*})}t)}
 \\&< \infty,
\end{align*}
for any $t\ge 1$.

 By It\^o's formula and $\textbf{(H7)}$, for any $\epsilon\in (0,\delta ^{*})$
 \begin{align*}
e^{\epsilon t}V(u(t;x))&=V(x)+\int_{0}^{t} \epsilon e^{\epsilon s}V(u(s))ds+\int_{0}^{t}e^{\epsilon s}[ _{B ^{\ast}} \langle A(u(s), \mathcal{L}_{u(s)}), \nabla V(u(s))\rangle _{B}
\\&~~~+\left \langle f(u_{s}, \mathcal{L}_{u_{s}}), \nabla V(u(s)) \right \rangle _{U_{1}}
+\frac{1}{2} \nabla ^{2}V(u(s))\left \| g(u_{s}, \mathcal{L}_{u_{s}}) \right \|^{2}_{\mathscr{L}(U_{2},U_{1})}]\text{d}s\\&~~~+\int_{0}^{t}e^{\epsilon s}  \left \langle \nabla V(u(s)),g(u_{s}, \mathcal{L}_{u_{s}})\text{d}W(s) \right \rangle _{U_{1}}
\\&\le V(x)-(\delta ^{*}-\epsilon)  \int_{0}^{t}e^{\epsilon s}V(u(s;x))\text{d}s+\epsilon^{-1}\delta_{0}e^{\epsilon t}+ \widetilde{\mathcal{M}}(t),
\end{align*}
where $ \widetilde{\mathcal{M}}(t) =\int_{0}^{t} e^{\epsilon s} \left \langle \nabla V(u(s)),g(u_{s}, \mathcal{L}_{u_{s}})\text{d}W(s) \right \rangle _{U_{1}}$. Then  by \eqref{a1} and Jensen's inequality, we obtain for any $t\ge 1$ and $c_{1}\in (0,\iota _{1})$
\begin{align}\label{a3}
 \mathbb{E} e^{c_{1}e^{-\epsilon t}\underset{z\in [0,t]}{\sup}\widetilde{\mathcal{M}}(t) }&\le \mathbb{E} e^{1+c_{1}e^{-\epsilon t}\widetilde{\mathcal{M}}(t) }\nonumber
 \\&\le e^{1+\frac{c_{1}^{2}\delta ^{2}_{2}}{\epsilon } }[\mathbb{E}e^{4c_{1}^{2}\delta ^{2}_{2} e^{-2\epsilon t}\int_{0}^{t} e^{2\epsilon s} V(u(s))\text{d}s}]^{\frac{1}{2}}\nonumber
 \\&\le e^{1+\frac{c_{1}^{2}\delta ^{2}_{2}}{\epsilon } }[\mathbb{E}e^{4c_{1}^{2}\delta ^{2}_{2} \frac{1-e^{-2\epsilon t}}{2\epsilon } \int_{0}^{t} V(u(s))\frac{2\epsilon}{1-e^{-2\epsilon t} } e^{-2\epsilon (t-s)} \text{d}s}]^{\frac{1}{2}}
 \\&\le e^{1+\frac{c_{1}^{2}\delta ^{2}_{2}}{\epsilon } }[\mathbb{E}\int_{0}^{t}e^{ \frac{2c_{1}^{2}\delta ^{2}_{2}}{\epsilon }  V(u(s))}\frac{2\epsilon}{1-e^{-2\epsilon t} } e^{-2\epsilon (t-s)} \text{d}s]^{\frac{1}{2}}\nonumber
 \\&\le e^{1+\frac{c_{1}^{2}\delta ^{2}_{2}}{\epsilon } }+\frac{2\epsilon e^{1+\frac{c_{1}^{2}\delta ^{2}_{2}}{\epsilon } }}{1-e^{-2\epsilon } }\mathbb{E}\int_{0}^{t}e^{ \frac{2c_{1}^{2}\delta ^{2}_{2}}{\epsilon }  V(u(s))} e^{-2\epsilon (t-s)} \text{d}s.\nonumber
\end{align}
Hence we have
\begin{align*}
\mathbb{E} (\sup_{z \in [0,t]}e^{c_{1}V(u(z;x) ) } )&=\mathbb{E} (e^{c_{1}\sup_{z \in [0,t]}V(u(z;x) ) } )
\\&\le e^{c_{1}V(x)+c_{1}\epsilon^{-1}\delta_{0}+\frac{c_{1}^{2}\delta ^{2}_{2}}{\epsilon } }+\frac{2\epsilon e^{c_{1}V(x)+c_{1}\epsilon^{-1}\delta_{0}+\frac{c_{1}^{2}\delta ^{2}_{2}}{\epsilon } }}{1-e^{-2\epsilon } }\mathbb{E}\int_{0}^{t}e^{ \frac{2c_{1}^{2}\delta ^{2}_{2}}{\epsilon }  V(u(s;x))} e^{-2\epsilon (t-s)} \text{d}s
\\&:=M_{x,\epsilon}+\frac{2\epsilon M_{x,\epsilon}}{1-e^{-2\epsilon } }\mathbb{E}\int_{0}^{t}e^{ \frac{2c_{1}^{2}\delta ^{2}_{2}}{\epsilon }  V(u(s))} e^{-2\epsilon (t-s)} \text{d}s
\\&\le M_{x,\epsilon}+\frac{\epsilon M^{2}_{x,\epsilon}}{(1-e^{-2\epsilon })^{2} }\int_{0}^{t}e^{-2\epsilon (t-s)} \text{d}s+\epsilon\mathbb{E}\int_{0}^{t}e^{ \frac{2c_{1}^{2}\delta ^{2}_{2}}{\epsilon }  V(u(s))} e^{-2\epsilon (t-s)} \text{d}s
\\&\le M_{x,\epsilon}+\epsilon\mathbb{E}\int_{0}^{t}\sup_{z \in [0,s]}e^{ \frac{2c_{1}^{2}\delta ^{2}_{2}}{\epsilon }  V(u(s))} e^{-2\epsilon (t-s)} \text{d}s.
\end{align*}
Let $ c_{1}=\frac{\epsilon}{2\delta ^{2}_{2} }$, then by  Gronwall inequality, there exist constants $c_{2}>0$ such that
\begin{align*}
\mathbb{E} (\sup_{z \in [0,t]}e^{c_{1}V(u(z;x) ) } )\le M_{x,\epsilon}+M_{x,\epsilon}\epsilon\mathbb{E}\int_{0}^{t}e^{-2\epsilon (t-s)}\text{d}s\le  e^{c_{2}(1+V( x)  )}.
\end{align*}
This completes the proof. $\Box $
~\\
\\\textbf{Remark 5.4.} It follows from \eqref{g10} that for any $\nu_{1},\nu_{2}\in\mathcal{P}  (U_{1} )$
\begin{align}\label{g25}
\begin{aligned}
\Pi _{V}(P^{*}_{t}\nu_{1} ,P^{*}_{t}\nu_{2}) \le e^{-\delta t } \int_{U_{1} }\int_{U_{1} }V( x -y  )\mu ^{*}(\text{d}y)\nu(\text{d}x ).
\end{aligned}
\end{align}
In addition, by \eqref{g15} and \eqref{g25},  we have
\begin{align}\label{g26}
\begin{aligned}
\left | P_{t}\Phi (x )- P_{t}\Phi (y )\right |&\le \left \| \Phi  \right \| _{LipV}\Pi _{V}( P^{*}_{t}\mu _{0}^{x },P^{*}_{t}\mu _{0}^{y })\\&\le M\left \| \Phi \right \| _{\mathcal{H}}e^{-\delta t }[\int _{U_{1} }\int _{U_{1}  }V(m-n)p(0,0,x ,\text{d}m )p(0,0, y,\text{d}n )]\\
&=M\left \| \Phi \right \| _{\mathcal{H}}e^{-\delta t }[\int _{\mathcal{R}  }V( x-n )p(0,0,y ,\text{d}n )]\\
&=M\left \| \Phi \right \| _{\mathcal{H}}e^{-\delta t }V( x-y).
\end{aligned}
\end{align}
Thus, by \eqref{g25} and \eqref{g26}, it is similar to  Lemma 4.2 in \cite{ref4},  and we can obtain that there exists a constant $M >0$ such that for any $\Phi\in C_{\mathcal{H} }(U_{1})$ and $x \in U_{1} $
\begin{align}\label{g27}
\begin{aligned}
 \left \| \Pi_{2} [\Phi] \right \| _{\mathcal{H} }\le M \left \| \Phi\right \| _{\mathcal{H} }^{2}.
\end{aligned}
\end{align}
\\\textbf{Lemma 5.5.}\emph{ Assume the conditions of Theorem 5.1 hold. Then for any $\Phi\in C_{\mathcal{H} }(U_{1} )$ with $(\Phi,\mu ^{*})=0$, we obtain
\begin{align*}
\begin{aligned}
 0\le \int _{U_{1}}\Pi_{2} [\Phi(x)]\mu ^{*}(\text{d}x ) =2 \int _{U_{1}}\Phi (x)\Pi_{1} [\Phi(x)]\mu ^{*}(\text{d}x )<\infty.
\end{aligned}
\end{align*}
\textbf{proof}} Similar to the proof of  Lemma 4.1 in \cite{ref4}, the above results can be obtained by \eqref{i1}, \eqref{g20}, \eqref{a5} and \eqref{a6}.

Subsequently, leveraging the aforementioned groundwork, we present the central limit theorem below.
~\\
\\\textbf{Theorem 5.6.}\emph{ Assume the conditions of Theorem 5.1 and \textbf{(H8)} with $\kappa =\frac{1}{2} $ hold.  For any $x\in U_{1}$ and $\Phi \in C_{\mathcal{H} }(U_{1} ) $ with $(\Phi,\mu ^{*})=0$, let $$\Lambda  =(\int _{U_{1}}\Pi_{2} [\Phi(x)]\mu ^{*}(\text{d}x ))^{\frac{1}{2} }\in [0,\infty ),$$ and we have the following conclusions:
\begin{enumerate}[\textbf{(i)}]
		\item When $\Lambda>0$, for $\varepsilon  \in [0,\frac{1}{4} )$, there exists an increasing function $\mathcal{I}  _{\varepsilon }:\mathbb{R}_{+}\times \mathbb{R} _{+}\to \mathbb{R} _{+}$ such that
\begin{align*}
\begin{aligned}
 \sup_{z\in \mathbb{R} } \left | \mathbb{P} (\frac{1}{\sqrt{t} }\int_{0}^{t}\Phi (u(s;x ))\text{d}s\le z  )-\Xi _{\Lambda } (z)\right | \le \mathcal{I}  _{\varepsilon }(\left \| \Phi \right \| _{\mathcal{H} },\left \| x  \right \|_{U_{1} } )t^{-\frac{1}{4}+\varepsilon  },
\end{aligned}
\end{align*}
for any $x \in U_{1} $ and $t\ge 1$;
\end{enumerate}
\begin{enumerate}[\textbf{(ii)}]
        \item When $\Lambda=0$, there exists an increasing function $\mathcal{I} :\mathbb{R} _{+}\times \mathbb{R} _{+}\to \mathbb{R} _{+}$ such that
\begin{align*}
\begin{aligned}
 \sup_{z\in \mathbb{R} }[ (\left | z \right |\wedge 1 ) \left | \mathbb{P} (\frac{1}{\sqrt{t} }\int_{0}^{t}\Phi(u(s;x ))\text{d}s\le z  )-\Xi _{0 } (z)\right | ]\le \mathcal{I}(\left \| \Phi \right \| _{\mathcal{H} },\left \| x  \right \|_{U_{1} } )t^{-\frac{1}{4} },
\end{aligned}
\end{align*}
for any $x \in U_{1} $ and $t\ge 1$,  where $$\Xi_{\Lambda }(r)=\frac{1}{\Lambda \sqrt{2\pi } } \int_{-\infty }^{r} e^{-\frac{s^{2}}{2\Lambda ^{2}} }\text{d}s,
\quad
\Xi _{0 }(r)=\left\{\begin{matrix}
  1,& r\ge0,\\
  0,& r<0.
\end{matrix}\right.$$
	\end{enumerate}
\textbf{Proof of (i):}} The  uniformly mixing of \eqref{r1} and Lemma 5.3 imply that the assumptions of Theorem 2.8 in \cite{ref5} hold, i.e.,
\begin{align*}
\begin{aligned}
 \mathbb{E} (\sup_{z \in [k,k+t]}e^{c_{1}V(u(z;x) ) } )\le e^{c_{2}(1+V( x)  )},
\end{aligned}
\end{align*}
for any $k\ge 0$ and $x\in U_{1}$. Thus, we can obtain that for $\bar{\Lambda } >0$ and $\varepsilon \in (0,\frac{1}{4} )$, there exists an increasing continuous function $\mathcal{K} _{\varepsilon }:\mathbb{R} _{+}\times \mathbb{R} _{+}\to \mathbb{R} _{+}$ such that for any $\Lambda\ge\bar{\Lambda }$ and $q>0 $
\begin{align}\label{g28}
\begin{aligned}
 &\sup_{z\in \mathbb{R} } \left | \mathbb{P} (\frac{1}{\sqrt{t} }\int_{0}^{t}\Phi(u(s;x ))\text{d}s\le z  )-\Xi _{\Lambda } (z)\right |\\ &\le t^{-\frac{1}{4}+\varepsilon  }\mathcal{K} _{\varepsilon }(\left \| \Phi \right \| _{\mathcal{H} },\left \| x  \right \|_{U_{1} } )+\Lambda^{-4q}\left [ t \right ]  ^{q(1-4\varepsilon)}\mathbb{E} \left | \frac{\Upsilon  ^{2}[ \mathcal{M} (\left [ t \right ])^{x}_{\Phi}]}{\left [ t \right ]  } -\Lambda ^{2} \right |^{2q}.
\end{aligned}
\end{align}
where $\left [ t \right ]$ is the integer part of $t$.
 Similar to \eqref{g15}, we have
\begin{align*}
\begin{aligned}
 \left | P_{t}\Pi_{2} [\Phi(x)] -\Lambda ^{2} \right | &=\left | P_{t}\Pi_{2} [\Phi(x)] -\int _{U_{1}}\Pi_{2} [\Phi(x)]\mu ^{*}(\text{d}x ) \right |\\
 &\le M\left \|\Pi_{2}[\Phi] \right \| _{\mathcal{H}}e^{-\delta t } [\frac{\delta_{0}}{c ^{*}\delta ^{*}(1-e^{- c ^{*}\delta ^{*}})}+V( cx   )],
\end{aligned}
\end{align*}
which implies that there exists a constant $M$ such that by   Lemma 2.1 of \cite{ref4} and \eqref{i1}, we have
\begin{align*}
\begin{aligned}
 \mathbb{E} \left | \frac{1}{n} \sum_{j=1}^{n}\|\Pi_{2}[\Phi(u(j;x) )]-\Lambda^{2} \right |^{2q} \le Mn^{-q}\left \| \Pi_{2}[\Phi] \right \|^{2q} _{\mathcal{H}}[\frac{\delta_{0}}{c ^{*}\delta ^{*}(1-e^{- c ^{*}\delta ^{*}})}+V( cx   )]^{2q}.
\end{aligned}
\end{align*}
Hence by \eqref{g23}, we arrive at
\begin{align}\label{g29}
\begin{aligned}
 \mathbb{E} \left | \frac{\Upsilon  ^{2}[\mathcal{M} (\left [ t \right ])^{x}_{\Phi}]}{\left [ t \right ]  } -\Lambda ^{2} \right |^{2q} \le M\left [ t \right ]^{-q}\left \| \Pi_{2}[\Phi] \right \|^{2q} _{\mathcal{H}}[\frac{\delta_{0}}{c ^{*}\delta ^{*}(1-e^{- c ^{*}\delta ^{*}})}+V( cx   )]^{2q}.
\end{aligned}
\end{align}
For an arbitrarily $\varepsilon\in (0, 1/4)$,
taking $q \ge  \frac{1}{16\varepsilon } $,  we have $\left[ t \right ]^{-4q\varepsilon } \le t^{-\frac{1}{4}+\varepsilon  }$ for any $t\ge 1$, which implies $$\mathcal{I} _{\varepsilon    }(\left \| \Phi \right \| _{\mathcal{H} },\left \| x  \right \|_{U_{1} } )=\mathcal{K} _{\varepsilon }(\left \| \Phi \right \| _{\mathcal{H} },\left \| x  \right \|_{U_{1} } )+M\Lambda^{-4q}\left \| \Phi \right \|^{4q} _{\mathcal{H}}[\frac{\delta_{0}}{c ^{*}\delta ^{*}(1-e^{- c ^{*}\delta ^{*}})}+V( cx   )]^{2q},$$
by \eqref{g28}, \eqref{g29} and Remark 5.4.
~\\
\\\textbf{\emph{Proof of (ii):}} When $\Lambda=0$, we can obtain  from Theorem 2.8 in \cite{ref5} that there exists an increasing continuous function $\mathcal{K}  :\mathbb{R} _{+}\times \mathbb{R} _{+}\to \mathbb{R} _{+}$ such that
\begin{align*}
\begin{aligned}
 &\sup_{z\in \mathbb{R} } \left | \mathbb{P} (\frac{1}{\sqrt{t} }\int_{0}^{t}\Phi(u(s;x ))\text{d}s\le z  )-\Xi _{0 } (z)\right |\\ &\le t^{-\frac{1}{4} } \mathcal{K} (\left \| \Phi \right \| _{\mathcal{H} },\left \| x  \right \|_{U_{1} })+\left [ t \right ] ^{-\frac{1}{2} } \left |\mathbb{E} \Upsilon  ^{2}[ \mathcal{M} (n)^{x}_{\Phi}]  \right |^{\frac{1}{2}}\\
 &\le t^{-\frac{1}{4} } \mathcal{K} (\left \| \Phi \right \| _{\mathcal{H} },\left \| x  \right \|_{U_{1} })+\left [ t \right ] ^{-\frac{1}{4} } M\left \| \Pi_{2}[\Phi] \right \|^{\frac{1}{2}} _{\mathcal{H}}[\frac{\delta_{0}}{c ^{*}\delta ^{*}(1-e^{- c ^{*}\delta ^{*}})}+V( cx   )]^{\frac{1}{2}}\\
 &\le\mathcal{I}(\left \| \Phi \right \| _{\mathcal{H} },\left \| x  \right \|_{U_{1} } )t^{-\frac{1}{4} },
\end{aligned}
\end{align*}
where $\mathcal{I}(\left \| \Phi \right \| _{\mathcal{H} },\left \| x  \right \|_{U_{1} } )=\mathcal{K} (\left \| \Phi \right \| _{\mathcal{H} },\left \| x  \right \|_{U_{1} } )+M \left \| \Phi \right \| _{\mathcal{H}}[\frac{\delta_{0}}{c ^{*}\delta ^{*}(1-e^{- c ^{*}\delta ^{*}})}+V( cx   )]^{\frac{1}{2}}$. The proof of Theorem 5.6 is complete. \quad $\Box$
\section*{Appendix \uppercase\expandafter{\romannumeral1}: The specific proof of  weak solutions(the step 1 of Theorem 3.2):}
In this section, we predominantly draw upon a segment of the proof provided in reference \cite{ref38}, Theorem 2.1. The specific details are outlined as follows:
~\\
\\\textbf{\emph{Proof:}} We first define the following  coordinate process:
\begin{align*}
x^{*}(t)\omega =\omega (t),
\end{align*}
where $\omega\in C(\mathbb{R}^{+} ,\mathbb{R}^{n} )$, and let $\mathcal{D} _{t}=\sigma \{\omega (s); 0\le s \le t\}$, hence $x^{*}(t)\omega$ is $\mathcal{D} _{t}$-adapted. Note that $x^{n}(t)$ with the initial condition $x^{n}(0) =x_{0}$ is the unique strong solution of  \eqref{r5}, which implies that
\begin{align*}
\mathcal{M} ^{n}(t)=x^{n}(t)-x_{0}-\int_{0}^{t} F^{n}(x^{n}(s),\mathcal{L}_{x^{n}(s)})\text{d}s
\end{align*}
is a martingale with the covariance  given by
\begin{align*}
\sum_{i=1}^{m} \int_{0}^{t} [G^{n}(x^{n}(s),\mathcal{L}_{x^{n}(s)})]_{ik}[G^{n}(x^{n}(s),\mathcal{L}_{x^{n}(s)})]_{jk}\text{d}s,
\end{align*}
$1\le i,j\le n$. Next let
\begin{align*}
\mathcal{M} ^{*,n}(t)=x^{*}(t)-x_{0}-\int_{0}^{t} F^{n}(x^{*}(s),\mathcal{L}_{x^{*}(s)})\text{d}s,
\end{align*}
then $\mathcal{M} ^{*,n}(t)$ is a martingale relative to $(\mathcal{L}_{x^{n}(t)},\mathcal{D} _{t})$ with the covariance
\begin{align*}
\left \langle \mathcal{M} ^{*,n}_{i},\mathcal{M} ^{*,n}_{j} \right \rangle(t) =\sum_{i=1}^{m} \int_{0}^{t} [G^{n}(x^{*}(s),\mathcal{L}_{x^{*}(s)})]_{ik}[G^{n}(x^{*}(s),\mathcal{L}_{x^{*}(s)})]_{jk}\text{d}s.
\end{align*}
Further by the property (\textbf{a}) and \eqref{r6}, let $n\to \infty $ and we obtain
\begin{align*}
\mathcal{M} ^{*}(t)=x^{*}(t)-x_{0}-\int_{0}^{t} F(x^{*}(s),\mathcal{L}_{x^{*}(s)})\text{d}s.
\end{align*}
Then for any $t>s$ and $\Gamma \in \mathcal{D}_{t} $, by  Problem 2.4.12 of \cite{ref66}, we have
\begin{align*}
\int _{\mathbb{R}^{n} }\mathcal{X}_{\Gamma }\mathcal{M} ^{*}(t)\text{d}\mathcal{L} _{x^{*}(t)}&=\lim_{n \to \infty} \int _{\mathbb{R}^{n} }\mathcal{X}_{\Gamma }\mathcal{M} ^{*,n}(t)\text{d}\mathcal{L} _{x^{*,n}(t)}\\&=\lim_{n \to \infty} \int _{\mathbb{R}^{n} }\mathcal{X}_{\Gamma }\mathcal{M} ^{*,n}(s)\text{d}\mathcal{L} _{x^{*,n}(s)}\\&=\int _{\mathbb{R}^{n} }\mathcal{X}_{\Gamma }\mathcal{M} ^{*}(s)\text{d}\mathcal{L} _{x^{*}(s)},
\end{align*}
where $\mathcal{X}_{\Gamma }$ represents the indicator function of $\Gamma$. This implies that $\mathcal{M} ^{*}(t)$ is a $\mathcal{L} _{x^{*}(s)}$-martingale. In addition, by the property (\textbf{a}) and \eqref{r6}, we get
\begin{align*}
\left \langle \mathcal{M} ^{*}_{i},\mathcal{M} ^{*}_{j} \right \rangle(t) =\sum_{i=1}^{m} \int_{0}^{t} [G(x^{*}(s),\mathcal{L}_{x^{*}(s)})]_{ik}[G(x^{*}(s),\mathcal{L}_{x^{*}(s)})]_{jk}\text{d}s
\end{align*}
for $1\le i,j\le n$. Based on Theorem II.7.1' of \cite{ref67}, it can be deduced that there exists an $m$-dimensional Brownian motion $B^{*}(t)$ on an extended probability space $(C(\mathbb{R}^{+} ,\mathbb{R}^{n} ),\mathcal{D} _{t},\mathcal{L}_{x^{*}(s)}))$ such that
\begin{align*}
\mathcal{M} ^{*}(t)=\int_{0}^{t} G(x^{*}(s),\mathcal{L}_{x^{*}(s)})\text{d}B^{*}(s),
\end{align*}
i.e.,
\begin{align*}
x^{*}(t)=x_{0}+\int_{0}^{t} F(x^{*}(s),\mathcal{L}_{x^{*}(s)})\text{d}s+\int_{0}^{t} G(x^{*}(s),\mathcal{L}_{x^{*}(s)})\text{d}B^{*}(s),
\end{align*}
hence $x^{*}(t)$ is  a weak solution to \eqref{r3}. The proof is complete.\quad $\Box$
\section*{Acknowledgments}
 The first author (S. Lu)  supported by Graduate Innovation Fund of Jilin University. The second author (X. Yang) was supported by  National Natural Science Foundation of China (12071175, 12371191). The third author (Y. Li) was supported by  National Natural Science Foundation of China (12071175 and 12471183).
%The authors are grateful to the referees and editor for their valuable suggestions and comments on this paper.
\section*{Data availability}
No data was used for the research described in the article

\section*{References}
\bibliographystyle{plain}
\bibliography{ref}

\end{document}